\documentclass[pre,twocolumn,twoside,byrevtex,superscriptaddress]{revtex4}

\usepackage{currfile}
\usepackage{url}
\usepackage{nameref}
\usepackage{hyperref}
\usepackage{graphicx,epsfig,verbatim,enumerate}
\usepackage{amssymb,amsmath}

\usepackage{xcolor}
\usepackage{color}

\definecolor{lightgrey}{rgb}{0.5,0.5,0.5}

\usepackage{rotating}
\usepackage{longtable}
\usepackage{array}

\usepackage{listings,color,setspace}

\newcommand{\sindex}[1]{}
\newcommand{\nindex}[1]{}

\newcommand{\etal}{\textit{et al.}}
\newcommand{\www}[1]{\url{#1}}

\newcommand{\pdiff}[2]{\frac{\partial #1}{\partial #2}}
\newcommand{\partialdiff}[2]{\frac{\partial #1}{\partial #2}}

\newcommand{\diff}[2]{\frac{{\rm d}#1}{{\rm d}#2}}

\newcommand{\rhoref}{\rho_{\text{ref}}}
\newcommand{\dphi}{\text{d}\phi}

\newcommand{\mbe}{\mathbf{\epsilon}}
\newcommand{\mbx}{\mathbf{x}}
\newcommand{\mby}{\mathbf{y}}
\newcommand{\mbd}{\mathbf{d}}
\newcommand{\mbB}{\mathbf{B}}
\newcommand{\mbW}{\mathbf{W}}
\newcommand{\mbR}{\mathbf{R}}
\newcommand{\mbH}{\mathbf{H}}
\newcommand{\mbK}{\mathbf{K}}
\newcommand{\mbP}{\mathbf{P}}
\newcommand{\mbZ}{\mathbf{Z}}
\newcommand{\mbw}{\mathbf{w}}
\newcommand{\mbX}{\mathbf{X}}
\newcommand{\mbY}{\mathbf{Y}}

\newcommand{\PreserveBackslash}[1]{\let\temp=\\#1\let\\=\temp}
\newcommand{\PBS}[1]{\let\temp=\\#1\let\\=\temp}

\newcommand{\plainlatexonly}[1]{}

\begin{document}

\title{
  Predicting Flow Reversals in a Computational Fluid Dynamics Simulated Thermosyphon using Data Assimilation

}

\author{Andrew J. Reagan}
\affiliation{Department of Mathematics \& Statistics, Vermont Complex Systems Center, Computational Story Lab, \& the Vermont Advanced Computing Core, The University of Vermont, Burlington, VT 05405}
\author{Yves Dubief}
\affiliation{School of Engineering, Vermont Complex Systems Center \& the Vermont Advanced Computing Core, The University of Vermont, Burlington, VT 05405}
\author{Peter Sheridan Dodds}
\affiliation{Department of Mathematics \& Statistics, Vermont Complex Systems Center, Computational Story Lab, \& the Vermont Advanced Computing Core, The University of Vermont, Burlington, VT 05405}
\author{Christopher M. Danforth}
\affiliation{Department of Mathematics \& Statistics, Vermont Complex Systems Center, Computational Story Lab, \& the Vermont Advanced Computing Core, The University of Vermont, Burlington, VT 05405}

\date{\today}

\begin{abstract}
A thermal convection loop is a annular chamber filled with water, heated on the bottom half and cooled on the top half.
With sufficiently large forcing of heat, the direction of fluid flow in the loop oscillates chaotically, dynamics analogous to the Earth's weather.
As is the case for state-of-the-art weather models, we only observe the statistics over a small region of state space, making prediction difficult.
To overcome this challenge, data assimilation (DA) methods, and specifically ensemble methods, use the computational model itself to estimate the uncertainty of the model to optimally combine these observations into an initial condition for predicting the future state.
Here, we build and verify four distinct DA methods, and then, we perform a twin model experiment with the computational fluid dynamics simulation of the loop using the Ensemble Transform Kalman Filter (ETKF) to assimilate observations and predict flow reversals.
We show that using adaptively shaped localized covariance outperforms static localized covariance with the ETKF, and allows for the use of less observations in predicting flow reversals.
We also show that a Dynamic Mode Decomposition (DMD) of the temperature and velocity fields recovers the low dimensional system underlying reversals, finding specific modes which together are predictive of reversal direction.
\end{abstract}

\maketitle

\section*{Introduction}

Prediction of the future state of complex systems is a fundamental challenge of science and engineering, and ultimately integral to the functioning of society.
Some of these systems include weather \cite{weather-violence2013}, health \cite{ginsberg2009a}, the economy \cite{sornette2006predictability}, marketing \cite{asur2010predicting} and transportation \cite{savely1972}.
For weather in particular, predictions are made using supercomputers integrating numerical weather models, projecting our current best guess of the atmospheric state into the future.
The accuracy of these predictions depends on the accuracy of the models themselves, and the quality of our knowledge of the current state of the atmosphere.

Model accuracy has improved with better meteorological understanding of weather processes and advances in computing technology \cite{bauer2015quiet}.
To solve the initial value problem, techniques developed over the past 50 years are now broadly known as {\em data assimilation} (DA).
Formally, data assimilation is the process of using all available information, including short-range model forecasts and physical observations, to estimate the current state of a system as accurately as possible \cite{yang2006}.
The best-guess of the current state is often referred to as the {\em analysis} state.

Here, we employ a fluid dynamics experiment as a test bed for improving numerical weather prediction algorithms, focusing specifically on data assimilation methods.
Our approach is inspired by the historical development of current methodologies, and provides a tractable system for rigorous analysis.
The experiment is a thermal convection loop, which by design simplifies our problem into the prediction of natural convection.
The thermosyphon, a type of natural convection loop or non-mechanical heat pump, can be likened to a toy model of climate \cite{harris2011predicting}.
The dynamics of thermal convection loops have been explored under both periodic \cite{keller1966} and chaotic \cite{welander1967,creveling1975stability,gorman1984,gorman1986,ehrhard1990dynamical,yuen1999,jiang2003,burroughs2005reduced,desrayaud2006numerical,yang2006,ridouane2010} regimes.
A full characterization of the computational behavior of a loop under flux boundary conditions by Louisos et. al. describes four regimes: chaotic convection with reversals, high Rayleigh number (Ra) aperiodic stable convection, steady stable convection, and conduction/quasi-conduction \cite{louisos2013}.
For the remainder of this work, we focus on the chaotic flow regime.

\section*{Physical Experiment and Computational Model}

The reduced order system describing a thermal convection loop was originally derived by Gorman \cite{gorman1986} and Ehrhard and M\"{u}ller \cite{ehrhard1990dynamical}.
Here we present this three dimensional system in non-dimensionalized form.
In Appendix B, we present a more complete derivation of these equations, following the derivation of Harris \cite{harris2011predicting}.
For the mean fluid velocity $\diff{x_1}{t}$, temperature difference between the 3 o'clock and 9 o'clock positions $\diff{x_2}{t}$ (also referred to presently as $\Delta T_{3-9}$), and deviation from conductive temperature profile $\diff{x_3}{t}$, these equations are:
\begin{align}
& \diff{x_1}{t} = \alpha (x_2 - x_1),\\
& \diff{x_2}{t} = \beta x_1 - x_2 (1 + Kh(|x_1|)) - x_1x_3,\\
  & \diff{x_3}{t} = x_1x_2 - x_3 (1 + Kh(|x_1|)) .\end{align}
The function $h(x)$ is a defined piece-wise analytic polynomial, and is provided in the full derivation as Equation \ref{eq:h_defined}.
The parameters $\alpha$, $\beta$, and $K$, along with scaling factors for time and each model variable can be fit to data using standard parameter estimation techniques.

\begin{figure}[h]
  \centering
  \includegraphics[width=0.5\textwidth]{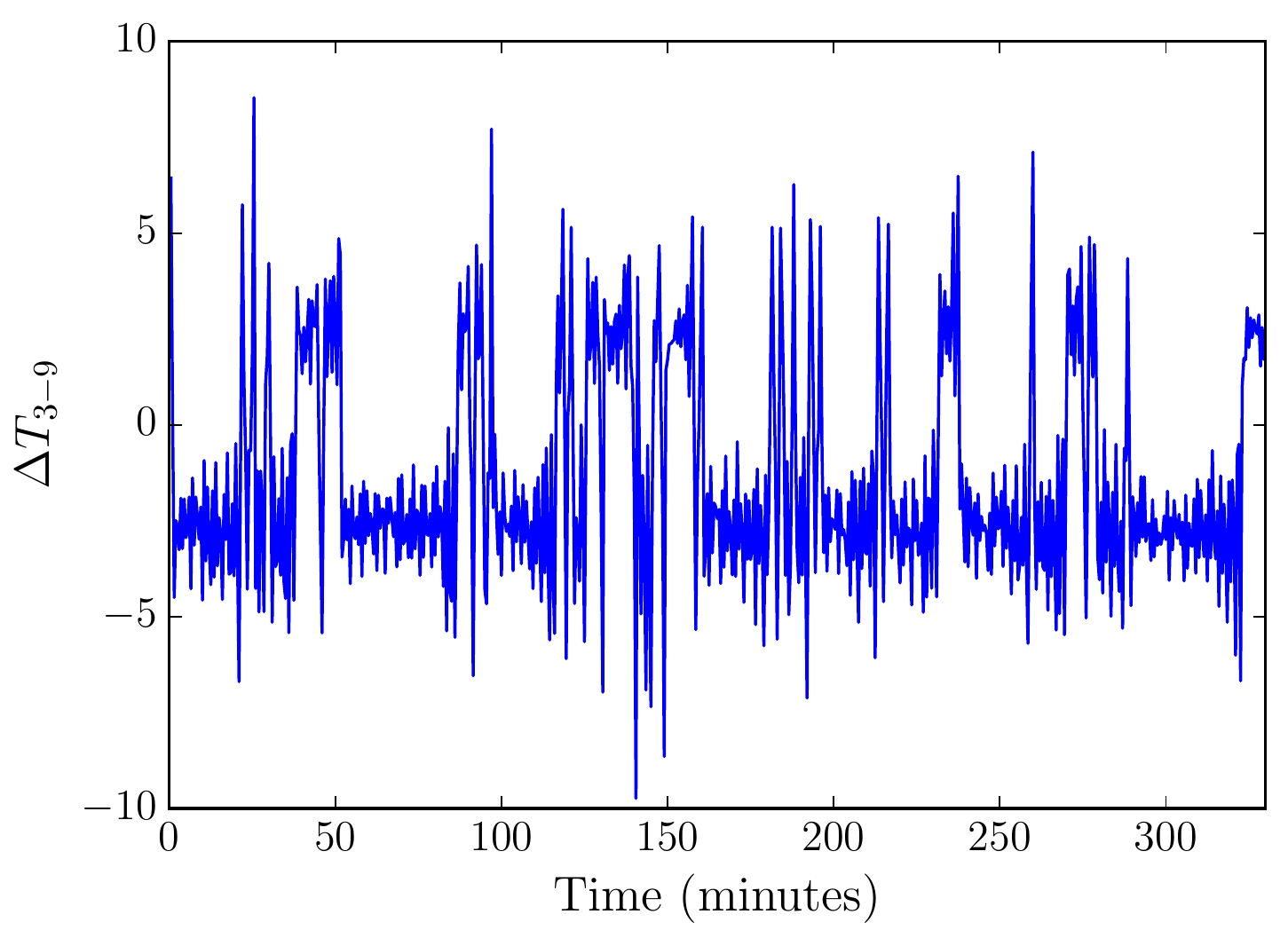}
  \caption[A time series of the physical thermosyphon, from the Undergraduate Honor's Thesis of Darcy Glenn {\protect \cite{glenn2013}}]{
    A time series of the physical thermosyphon, from the Undergraduate Honor's Thesis of Darcy Glenn {\protect \cite{glenn2013}}.
    The temperature difference (plotted) is taken as the difference between temperature sensors in the 3 and 9 o'clock positions.
    The sign of the temperature difference indicates the flow direction, where positive values are clockwise flow.
      }
  \label{fig:thermosyphon-physical-timeseries}
\end{figure}

\begin{figure}[h]
  \centering
  \includegraphics[width=0.45\textwidth]{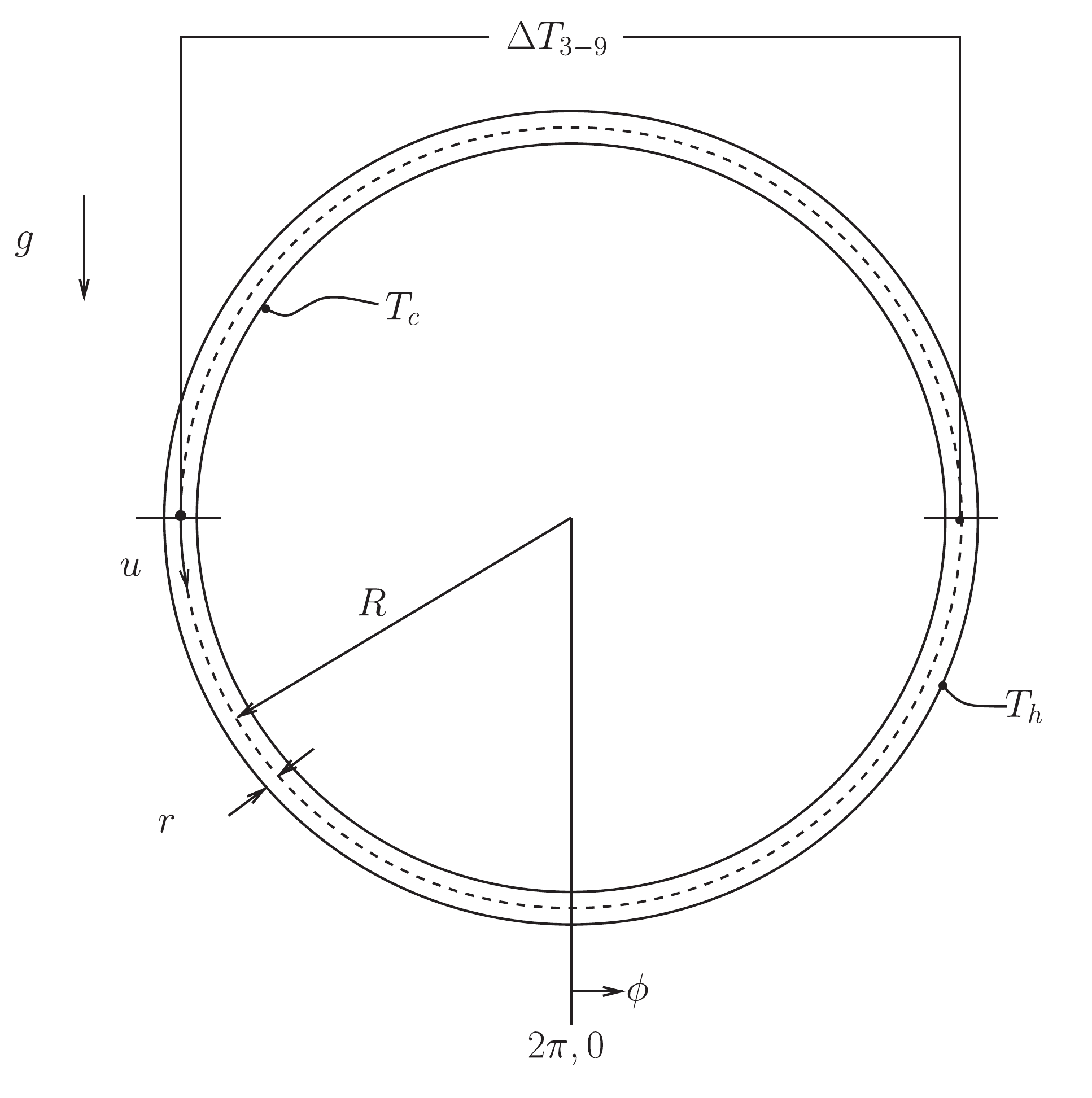}
  \caption[Schematic of the experimental, and computational, setup from Harris \etal (2012)]{
    Schematic of the experimental, and computational, setup from Harris \etal (2012).
    The loop radius is given by $R$ and inner radius by $r$.
    The top temperature is labeled $T_c$ and bottom temperature $T_h$, gravity $g$ is defined downward, the angle $\phi$ is prescribed from the 6 o'clock position, and temperature difference between 3 o'clock and 9 o'clock positions $\Delta T_{3-9}$ is labeled.
  }
  \label{fig:thermosyphon-schematic}
\end{figure}

Operated by Dave Hammond, UVM's Scientific Electronics Technician, the experimental thermosyphons access the chaotic regime of state space found in the principled governing equations.
We quote the detailed setup from Darcy Glenn's undergraduate thesis \cite{glenn2013} and provide Fig. \ref{fig:thermosyphon-schematic} for details of the experiment:
\begin{quote}
The [thermosyphon] is a bent semi-flexible plastic tube with a 10-foot heating rope wrapped around the bottom half of the upright circle.
The tubing used is light-transmitting clear THV from McMaster-Carr, with an inner diameter of 7/8 inch, a wall thickness of 1/16 inch, and a maximum operating temperature of 200F.
The outer diameter of the circular thermosyphon is 32.25 inches.
Together, the tubing inner diameter and outer diameter of the thermosyphon produce a ratio of approximately 1:36.
There are 1 inch 'windows' when the heating cable is coiled in a helix pattern around the outside of the tube, so the heating is not exactly uniform.
The bottom half is then insulated using aluminum foil, which allowed fluid in the bottom half to reach 176F.
A forcing of 57 V, or 105 Watts, is required for the heating cable so that chaotic motion is observed.
Temperature is measured at the 3 o'clock and 9 o'clock positions using unsheathed copper thermocouples from Omega.
\end{quote}
We confirm that the experiment accesses the chaotic regime of state space using a time series of the temperature difference as measured at the 3 o'clock and 9 o'clock positions in Fig. \ref{fig:thermosyphon-physical-timeseries}.
We first test our ability to predict this experimental thermosyphon using synthetic data.

We perform all computational simulations of the thermal convection loop with the open-source finite volume C++ library OpenFOAM \cite{jasak2007}.
The open-source nature of this software enables its integration with the data assimilation framework that our present work provides.

We consider the incompressible Navier-Stokes equations with the Boussinesq approximation to model the flow of water inside a thermal convection loop.
For brevity, we omit the equations themselves, and include them in the Appendix.
The solver in OpenFOAM that we use, with some modification, is \verb|buoyantBoussinesqPimpleFoam|.
Solving is accomplished by the Pressure-Implicit Split Operator (PISO) algorithm \cite{issa1986solution}.
We find that modification of the code is necessary for laminar operation.

\begin{figure}[h]
  \centering
  \includegraphics[width=0.45\textwidth]{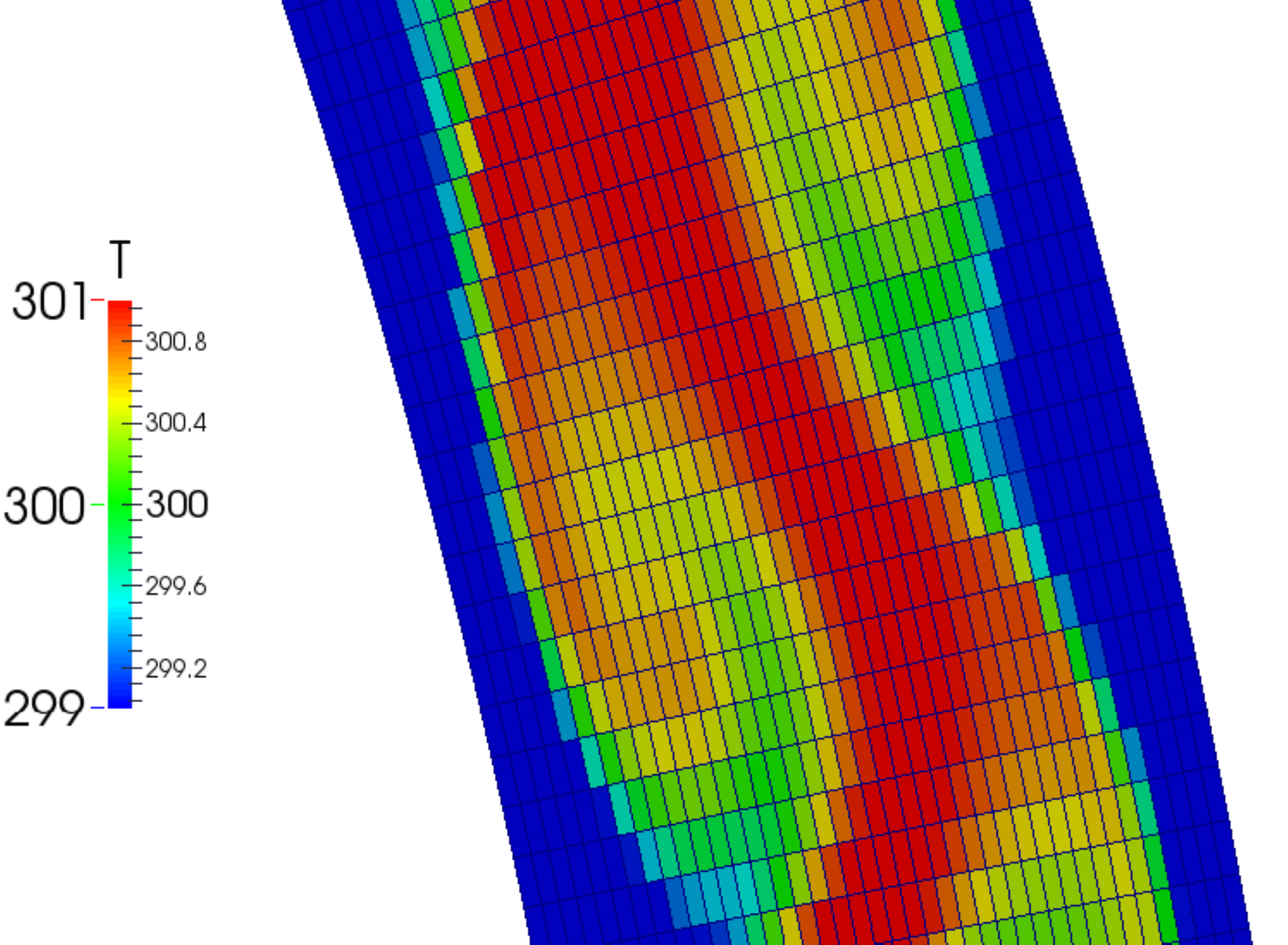}
  \caption[A snapshot of the mesh used for CFD simulations]{
    A snapshot of the mesh used for CFD simulations.
    Shown is an initial stage of heating for a fixed value boundary condition, 2D, laminar simulation with a mesh of 40000 cells without wall refinement with walls heated at 340K on the bottom half and cooled to 290K on the top half.
    The cells have been colored with a truncated temperature range (299--301K) to highlight the flow structures.
  }
  \label{fig:CFDmesh1}
\end{figure}

\begin{figure}[h]
  \centering
  \includegraphics[width=0.45\textwidth]{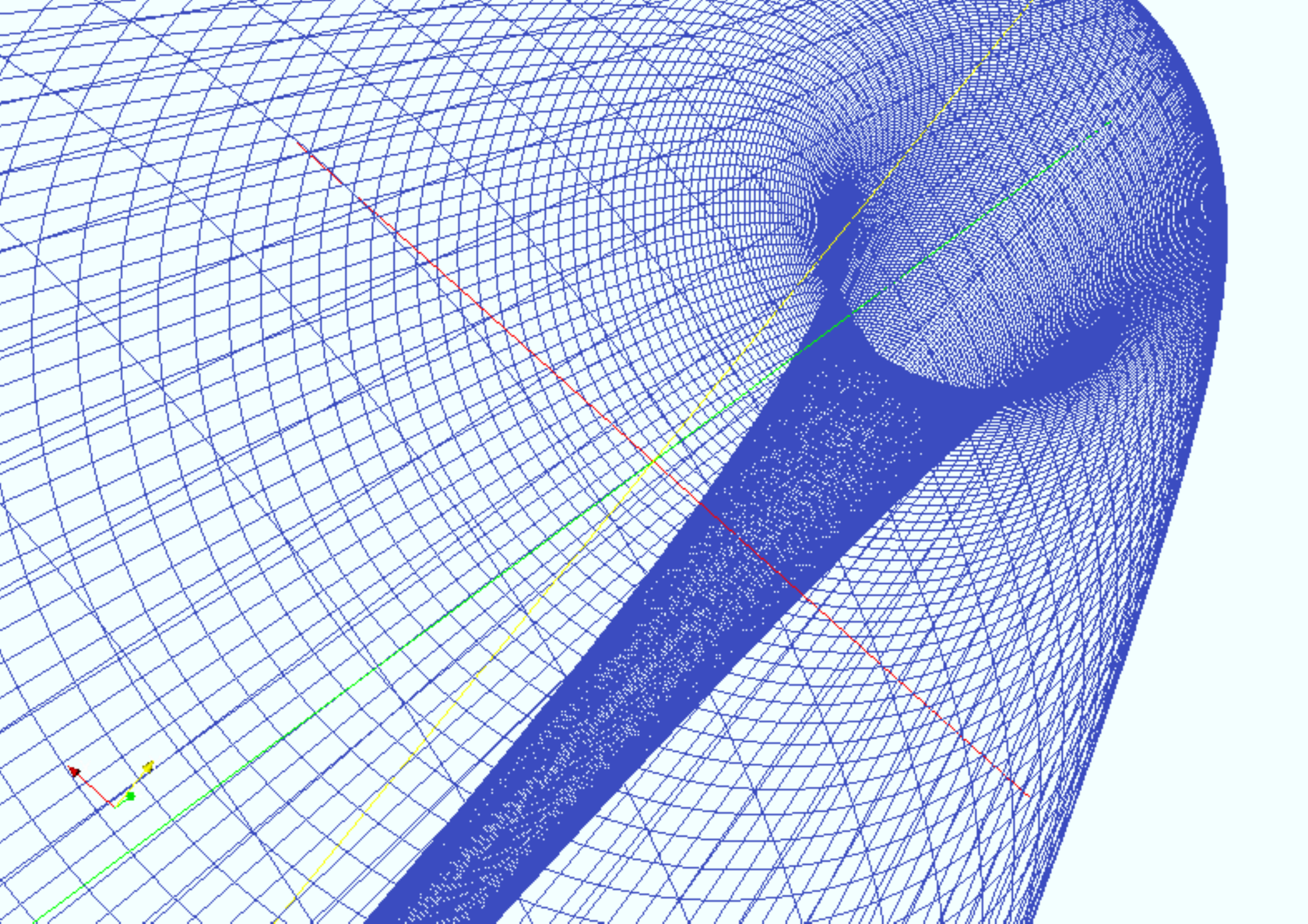}
  \caption[The 3D mesh viewed as a wire-frame from within]{
    The 3D mesh viewed as a wire-frame from within.
    Here there are 900 cells in each slice (not shown), for a total mesh size of 81,000 cells.
    Simulations using this computational mesh are prohibitively expensive for use in a real time ensemble forecasting system, but are possible offline.
  }
  \label{fig:CFDmesh2}
\end{figure}

We create both 2-dimensional and 3-dimensional meshes using OpenFOAM's native meshing utility \verb|blockMesh| shown in Figures \ref{fig:CFDmesh1} and \ref{fig:CFDmesh2}.
After creating a mesh, we refine the mesh near the walls to capture boundary layer phenomena and renumber the mesh for solving speed.
We use the \verb|refineWallMesh| utility to refine the mesh near walls, and the \verb|renumberMesh| utility to renumber the mesh.
The resulting 2D mesh contains 80,000 points (80 across the diameter and 1000 around).

Available boundary conditions (BCs) we find to be stable in OpenFOAM's solver are constant gradient, fixed value conditions, and turbulent heat flux.
Constant gradient simulations are stable, but the behavior is empirically different from our physical system.
While it is possible that a fixed value BC is acceptable due to the thermal diffusivity and thickness of the walls of the experimental setup, we find that this is also inadequate.
Simulations with a turbulent heat flux BC implemented through the \verb|externalWallHeatFluxTemperature| library are unstable with the laminar turbulence model we use and resulted in physically unrealistic results.
We employ the third-party library \verb|groovyBC| to use a gradient condition that computes the flux using a fixed external temperature $T_\text{inf}$ and fixed wall heat transfer coefficient $h$ as $$ - \partialdiff{T}{x_j} = h \left( T-T_\text{inf} \right)$$ where we choose $h$ to be the reference value for aluminum (the material used in the experimental setup).

With the mesh, BCs, and solver chosen, we now simulate the flow.
From the data of $T,\phi,\vec{u}$ and $p$ that are saved at each timestep (temperature, cell face flux, velocity, and pressure, respectively), we extract the mass flow rate and average temperature at the $12,3,6$ and $9$ o'clock positions on the loop.
Since $\phi$ is saved as a face-value flux, we compute the mass flow rate over the cells $i$ of top (12 o'clock) slice as
\begin{equation} \sum _i\phi_{f(i)} \cdot v_i \cdot \rho_i\end{equation}
where $f(i)$ corresponds the face perpendicular to the loop angle at cell $i$ and $\rho$ is reconstructed from the Boussinesq approximation $\rho = \rhoref (1-\beta(T-T_\text{ref}))$.

\section*{Methods}

\subsection*{Data Assimilation}

We perform initial tests of the data assimilation algorithms described here with the Lorenz '63 system, which is analogous to the above equations with Lorenz's $\beta = 1$, and $K = 0$.
The canonical choices of $\sigma = 10, \beta = 8/3$ and $\rho = 28$ produce the well known butterfly attractor, and we use these values for all examples here.
From these tests, we will find the optimal data assimilation parameters (inflation factors) for predicting time series with this system.
Having done so, we then focus our efforts on making prediction using computational fluid dynamics models.

We first implement the 3D-Var filter.
Simply put, 3D-Var is the variational (cost-function) approach to finding the analysis.
It has been shown that 3D-var solves the same statistical problem as optimal interpolation (OI) \cite{lorenc1986analysis}.
The usefulness of the variational approach comes from the computational efficiency, when solved with an iterative method.
Specifically, the multivariate 3D-Var amounts to finding the $\mbx _a$ that minimizes the cost function
\begin{equation} J(\mbx) = (\mbx - \mbx_b) ^T \mbB ^{-1} (\mbx - \mbx_b) + (\mby_o + H(\mbx))^T\mbR (\mby_o - H(\mbx)) .\end{equation}

Next, we implement the ``gold-standard'' Extended Kalman Filter (EKF).
The tangent linear model (TLM) is precisely the model (written as a matrix) that transforms a perturbation at time $t$ to a perturbation at time $t+\Delta t$, analytically equivalent to the Jacobian of the model.
Using the notation of Kalnay \cite{kalnay2003}, this amounts to making a forecast with the nonlinear model $M$, and updating the error covariance matrix $\mbP$ with the TLM $L$, and adjoint model $L^T$:

\begin{align*} \mbx^f (t_i) &= M _{i-1} [\mbx ^a (t_{i-1} ) ],\\
\mbP^f (t_i ) &= L_{i-1} \mbP^a (t_{i-1} ) L^T _{i-1} + \mathbf{Q} (t_{i-1} ) \end{align*}

where $\mathbf{Q}$ is the noise covariance matrix (model error).
In the experiments with Lorenz '63 presented in this section, $\mathbf{Q} = 0$ since our model is perfect.
In numerical weather prediction, $\mathbf{Q}$ must be approximated, e.g., using statistical moments on the analysis increments \cite{danforth2007estimating,li2009accounting,danforth2008using}.

The analysis step is then written as (for $H$ the observation operator):
\begin{align} \mbx^a (t_i ) &= \mbx^f (t_i) + \mbK_i \mbd_i,\\
\mbP^a (t_i) &= (\mathbf{I} - \mbK_i \mbH_i )\mbP^f (t_i) \end{align}
where
\[ \mbd_i = \mby_i^o - \mbH[x^f (t_i) ] \]
is the innovation. We compute the Kalman gain matrix to minimize the analysis error covariance $P^a _i$ as
\[ \mbK_i = \mbP^f (t_i) \mbH_i ^T [ \mbR_i + \mbH_i \mbP^f (t_i) \mbH^T ] ^{-1} \]
where $\mbR_i$ is the observation error covariance.
Since we are making observations of the truth with random normal errors of standard deviation $\mbe$, the observational error covariance matrix $\mbR$ is a diagonal matrix with $\epsilon$ along the diagonal.
The most difficult (and most computationally expensive) part of the EKF is deriving and integrating the TLM.
For this reason, the EKF is not used operationally, and later we will turn to statistical approximations of the EKF using ensembles of model forecasts.
With our CFD model we have no such TLM, and we provide more detail on the TLM approaches applicable to the Lorenz '63 system in Appendix C.

The computational cost of the EKF is mitigated through the approximation of the error covariance matrix $\mbP_f$ from the model itself, without the use of a TLM.
One such approach is the use of a forecast ensemble, where a collection of models (ensemble members) are used to statistically sample model error propagation.
With ensemble members spanning the model analysis error space, the forecasts of these ensemble members are then used to estimate the model forecast error covariance.

The only difference between this approach and the EKF, in general, is that the forecast error covariance $\mbP^f$ is computed from the ensemble members, without the need for a tangent linear model:
\[ \mbP^f \approx \frac{1}{K-2} \sum _{k\neq l} \left ( \mbx_k ^f - \overline{\mbx} ^f _l \right ) \left (\mbx_k ^f - \overline{\mbx} ^f _l \right ) ^T .\]

The ETKF introduced by Bishop is one type of square root filter, and we present it here to provide background for the formulation of the LETKF \cite{bishop2001adaptive}.
For a square root filter in general, we begin by writing the covariance matrices as the product of their matrix square roots.
Because $\mbP_a$ and $\mbP_f$ are symmetric positive-definite (by definition), we can write
\begin{equation} \mbP_a = \mbZ_a \mbZ_a^T ~~,~~~ \mbP_f = \mbZ_f \mbZ_f^T \end{equation}

where $\mbZ_a$ and $\mbZ_f$ are the matrix square roots of $\mbP_a$ and $\mbP_f$.
We are not concerned that this decomposition is not unique, and note that $\mbZ$ must have the same rank as $\mbP$ which will prove computationally advantageous.
The power of the SRF is now seen as we represent the columns of the matrix $\mbZ_f$ as the difference from the ensemble members from the ensemble mean, to avoid forming the full forecast covariance matrix $\mbP_f$.
The ensemble members are updated by applying the model $M$ to the states $\mbZ_f$ such that an update is performed by
\begin{equation} \mbZ_f = M \mbZ_a .\end{equation}
To summarize, the steps for the ETKF are to (1) form $\mbZ_f^T\mbH^T\mbR^{-1}\mbH\mbZ_f$, assuming that computing $\mathbf{R}^{-1}$ is easy, and (2) compute its eigenvalue decomposition, and apply it to $\mbZ_f$.

The LEKF implements a strategy that becomes important for large simulations: localization.
Namely, the analysis is computed for each grid-point using only local observations, without the need to build matrices that represent the entire analysis space.
Localization removes long-distance correlations from $\mathbf{B}$ and allows greater flexibility in the global analysis by allowing different linear combinations of ensemble members at different spatial locations \cite{kalnay2007a}.
The general formulation of the LEKF by Ott goes as follows, quoting directly from \cite{ott2004local}:
\begin{enumerate}
\item Globally advance each ensemble member to the next analysis timestep. Steps 2--5 are performed for each grid point.
\item Create local vectors from each ensemble member.
\item Project that point's local vectors from each ensemble member into a low dimensional subspace as represented by perturbations from the mean.
\item Perform the data assimilation step to obtain a local analysis mean and covariance.
\item Generate local analysis ensemble of states.
\item Form a new global analysis ensemble from all of the local analyses.
\item Wash, rinse, and repeat.
\end{enumerate}

Proposed by Hunt \etal (2007) with the stated objective of computational efficiency, the LETKF is named from its most similar algorithms from which it draws \cite{hunt2007efficient}.
With the formulation of the LEKF and the ETKF given, the LETKF can be  described as a synthesis of the advantages of both of these approaches.
The LETKF is the method sufficiently efficient for implementation on the full OpenFOAM CFD model of 240,000 model variables, and so we present it in more detail and follow the notation of Hunt \etal (2007). 
As in the LEKF, we explicitly perform the analysis for each grid point of the model.
The choice of observations to use for each grid point can be selected a priori, and tuned adaptively.
Starting with a collection of background forecast vectors $\{ \mbx_{b(i)}:\,i=1,\ldots,k \}$, we perform steps 1 and 2 in a global variable space, then steps 3--8 for each grid point:
\begin{enumerate}
\item Apply $H$ to $\mbx_{b(i)}$ to form $\mby_{b(i)}$, average the $\mby_b$ for $\overline{\mby_b}$, and form $\mbY b$.
\item Similarly form $\mbX_b$. Now for each grid point:
\item Form the local vectors.
\item Compute $\mathbf{C}=(\mbY_b)^T\mbR^{-1}$ (perhaps by solving $\mbR \mathbf{C}^T = \mbY_b$.
\item Compute $\tilde{\mbP}_a = \left( (k-1)\mathbf{I} / \rho + \mathbf{C} \mbY _b \right ) ^{-1}$ where $\rho > 1$ is a tun-able covariance inflation factor.
\item Compute $\mbW_a = \left ( (k-1) \tilde{\mbP} _a \right ) ^{1/2}$.
\item Compute $\overline{\mbw} _a  = \tilde{\mbP}_a \mathbf{C} \left ( \mby_o - \overline{\mby} _b \right )$ and add it to the column of $\mbW_a$.
\item Multiply $\mbX_b$ by each $\mbw_{a(i)}$ and add $\tilde{\mbx}_b$ to get $\left\{ \mbx_{a(i)}:\,i=1,\ldots,k\right \}$ to complete each grid point.
\item Combine all of the local analysis into the global analysis.
\end{enumerate}
We implement the LETKF on our mesh using the full 80 cells across with zone sizes of center 10, and sides 15, resulting in 3200 local variables for 100 zones.
In parallel, these 100 local computations can all be carried out simultaneously over an arbitrary number of processors.

\subsection*{Adaptive covariance localization}

\begin{figure}[h]
  \centering
  \includegraphics[width=0.45\textwidth]{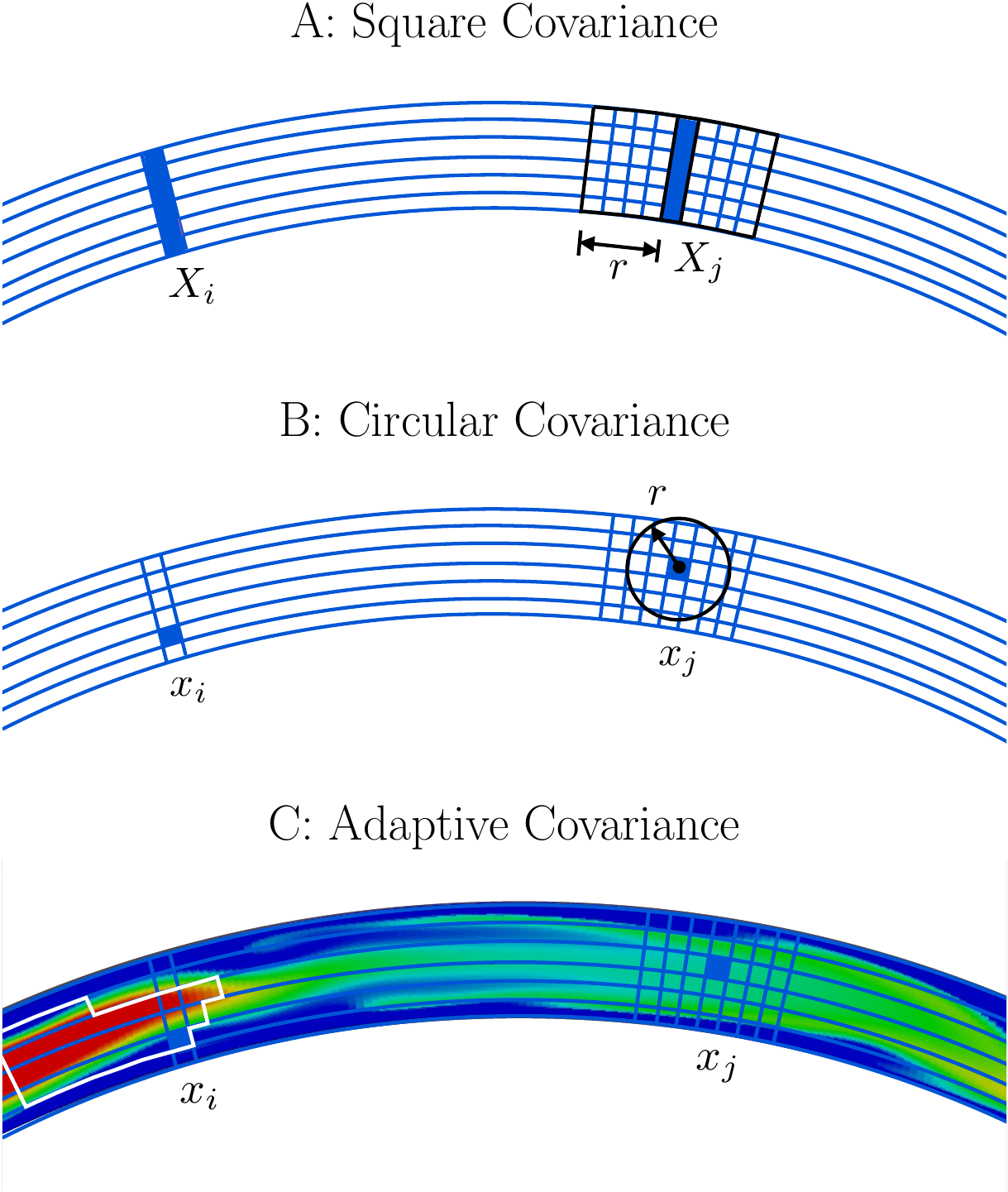}
  \caption[]{
    Schematic of the adaptive covariance localization.
    In Panel A we see a zonal (square) covariance that is most similar to the covariance used for both control experiments and sliding covariance experiments.
    Panel B shows a localized covariance using a ``local radius'', and Panel C shows an idealized, fully adaptive covariance.
    While we are motivated by localization around flow structures like Panel C, we simply shift the covariance in Panel A so that our method is most general and computationally efficient.
  }
  \label{fig:covariance-localization-schematic}
\end{figure}

Using the ``square'' sections of the loop to localize, we shift the zone to the left or right to follow the dominate flow direction at the center of that local window.
In Fig. \ref{fig:covariance-localization-schematic} a schematic of localization using square, circular, and adaptive location shows a situation in which adaptive localization will potentially capture more relevant information for finding the analysis state of any given cell.
As we note in the caption of Fig. \ref{fig:covariance-localization-schematic}, while we are motivated by localization around flow structures like Panel C, we simply shift the covariance in Panel A so that our method is most general and computationally efficient.

Denote the velocity vector of cells on a perpendicular slice of the loop at $\vec{U}$, the tangent vector to the slice $\vec{U}$ by $\vec{T}$, the zone width as $z_{\text{max}}$ and then the localization shift $\alpha_{\text{local}}$ for that slice of the loop is taken to be
\begin{equation} \alpha_{\text{local}} = \text{floor} \left( (\vec{U} \cdot \vec{T})/\max (U) \times z_{\text{max}} \right) . \end{equation}

\subsection*{Dynamic mode decomposition}

We employ the ``standard'' algorithm of Tu to compute the Dynamic Mode Decomposition \cite{tu2013dynamic}.
Tu's ``standard'' algorithm is as follows with $X$ and $Y$ taken as the first and last $N-1$ columns of the snapshot matrix $D$:
\begin{align*} X &= U\Sigma V \tag*{(Take SVD of $X$.)}\\
  \tilde{A} &= U^T Y V \Sigma ^{-1} \tag*{(Build the $A$ matrix.)}\\
  \tilde{A}w &= \lambda w \tag*{(Compute eigenvectors and values.)}\\
  \hat{\theta}w &= U w \tag*{(Compute corresponding modes.)}\end{align*}

Given a system state $U^*$ we project this state onto the DMD basis by taking the real part of $\Phi = \text{re}\left (U^*\cdot w\right)$ and use the psuedoinverse to compute the projection as $$(\Phi^T \cdot \Phi)^{-1} \cdot \Phi ^T \cdot U^*.$$
This projection is a vector which contains the linear coefficients on the basis of DMD modes for the given state.

\section*{Results}

\subsection*{Data assimilation}
\label{data_assimilation_section}

\begin{figure}[h]
  \centering
  \includegraphics[width=0.50\textwidth]{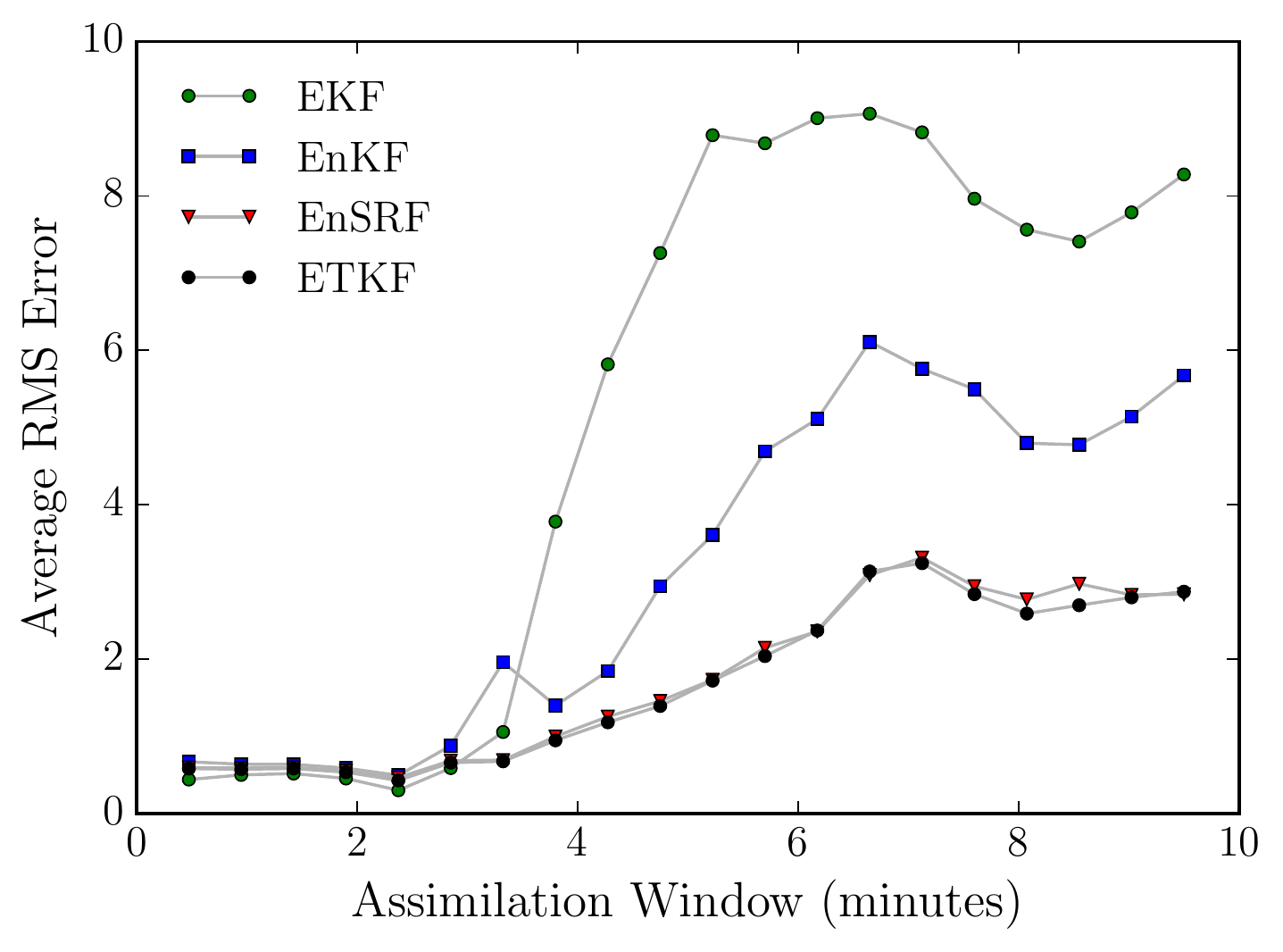}
  \caption[The RMS error is reported for our EKF and EnKF filters]{
    The RMS error (not scaled by climatology) for our EKF and EnKF filters.
    Error is measured as the difference between forecast and truth at the end of an assimiliation window for the latter 2500 assimiliation windows in a 3000 assimilation window Lorenz '63 run.
    Error is measured in the only observed variable, $x_1$.
    Increasing the assimilation window led to an decrease in predictive skill, as expected.
  }
  \label{fig:window_test}
\end{figure}

We confirm the performance of the DA methods described above by testing each (on the Lorenz '63 system) for increasingly long times between observations, by increasing the DA window length in Fig. \ref{fig:window_test}.
As the time between observations increases, the nonlinearity of the Lorenz '63 system results in the failure of the EKF and difficultly for the EnKF with small ensemble size.
The ETKF and EnSRF perform the best of the methods tested and we chose the ETKF for future use with the CFD model.

\begin{figure}[h]
  \centering
  \includegraphics[width=0.48\textwidth]{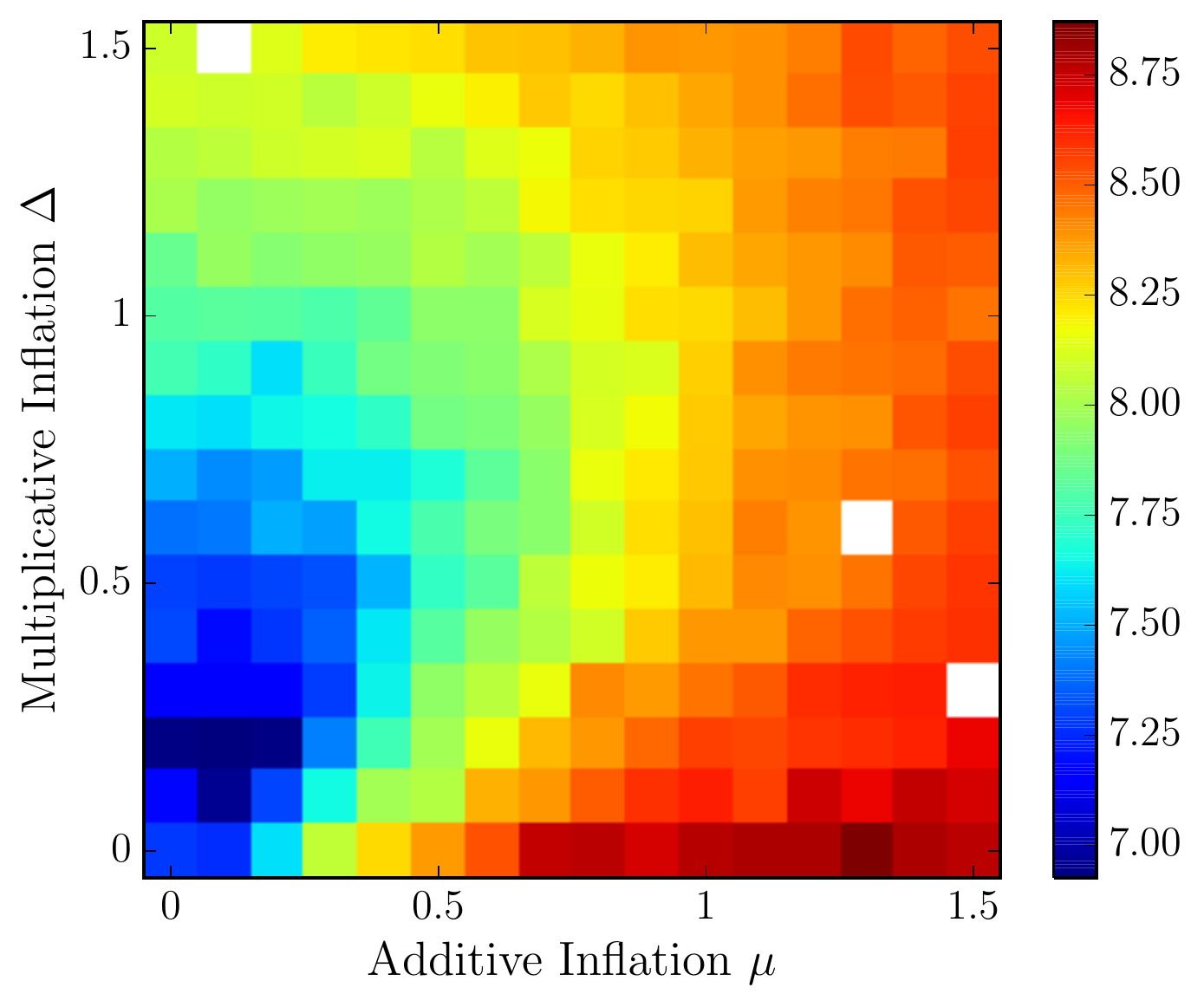}
  \caption[The RMS error averaged over 100 model runs of length 1000 windows is reported for the ETKF for varying additive and multiplicative inflation factors]{
    The RMS error averaged over 100 model runs of length 1000 windows is reported for the ETKF for varying additive and multiplicative inflation factors $\Delta$ and $\mu$.
    Each of the 100 model runs starts with a random IC, and the analysis forecast starts randomly.
    The window length here is 390 seconds.
    The filter performance RMS is computed as the RMS value of the difference between forecast and truth at the assimilation window for the latter 500 windows, allowing a spin-up of 500 windows.
  }
  \label{fig:ETKF_cov_tuning_390s}
\end{figure}

The results in Fig. \ref{fig:window_test} rely on tuned covariance inflation, both additive and multiplicative, pre-computed for each window and DA technique.
We choose optimal additive inflation $\mu$ and multiplicative inflation $\Delta$ by selecting for the lowest error in an exhaustive search through a maximum factors of $1.5$ in each, an example is shown in Fig. \ref{fig:ETKF_cov_tuning_390s}.
We use these optimal data assimilation parameters (inflation factors) for the remainder of this work.

\subsection*{Limited observations \& adaptive covariance}

An initial test of prediction skill with limited observations in a twin model experiment showed that we needed 1000 spatial measurements of the temperature to predict flow reversals within 1 assimilation window.
In an attempt to decrease the required observations to a experimentally realizable number, we implement a simple, adaptively localized covariance for data assimilation.
Since we first saw a modest improvement in the prediction skill with full temperature observations, we hope that this improvement increases and is sufficient to get down to needing as few as 32 observations to predict reversals 1 assimilation window (of length 10 seconds) into the future.

\begin{figure}[h]
  \centering
  \includegraphics[width=0.48\textwidth]{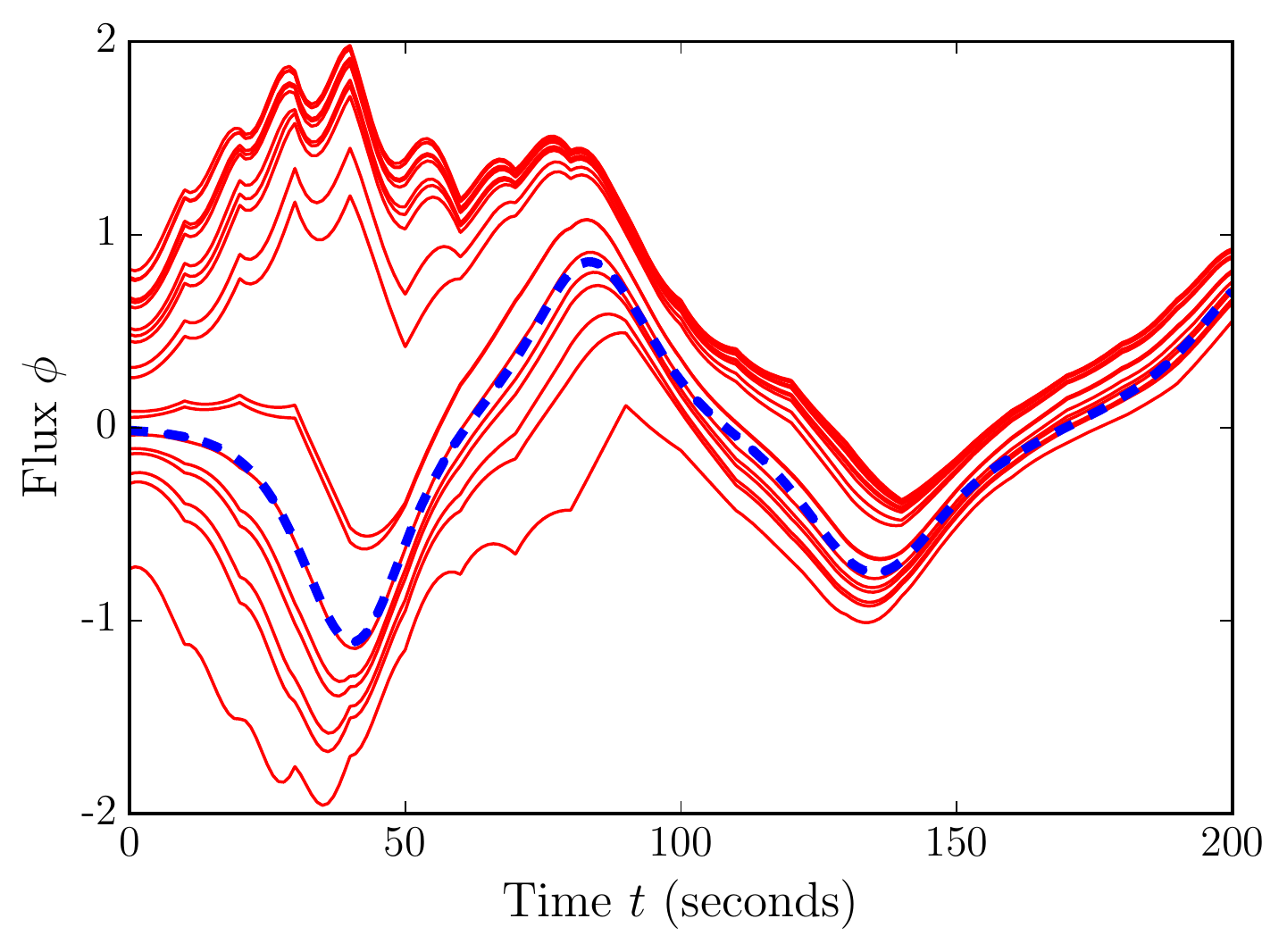}
  \caption[]{
    Convergence of 20 ensembles using sliding windows, starting from initially random states.
    Here, as in most of the experiments, only temperature is observed and assimilated.
    Flux is computed as in Equation 4, on the left hand side of the thermosyphon, and scaled by a factor of $10^8$.
    Assimilation takes place every 10 model seconds.
  }
  \label{fig:sliding_spag}
\end{figure}

\begin{figure}[h]
  \centering
  \includegraphics[width=0.45\textwidth]{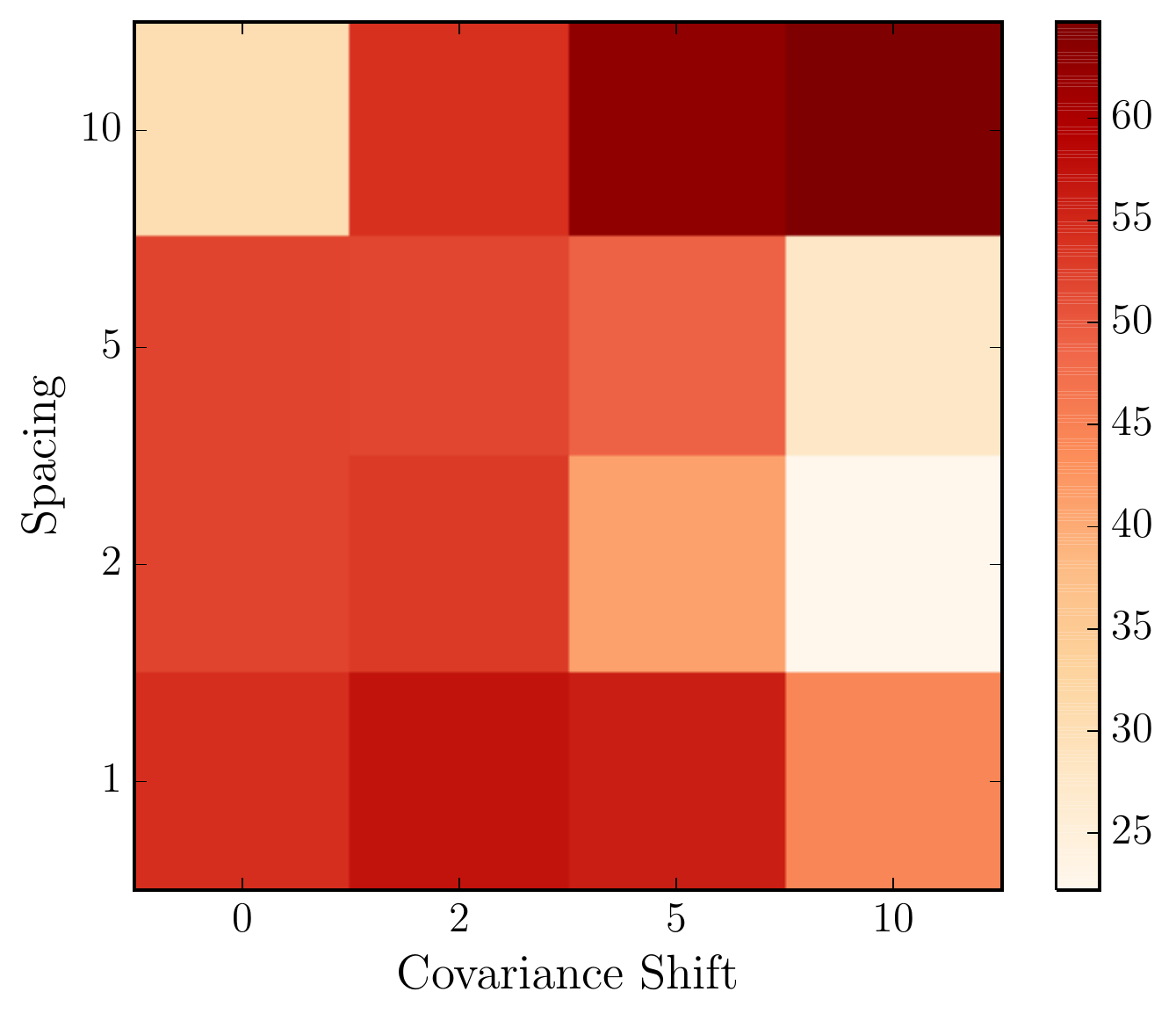}
  \caption[]{
    Prediction skill as fraction of reversals that we correctly predicted across different numbers of observations and sliding windows of localized covariance.
    Decreasing observation density makes the prediction problem more difficult while at the same time make the data assimilation stable numerically, and we see a decrease in prediction skill with no covariance shifting.
    With covariance shifting, skill improves for each observational density and most dramatically with less observation density.
  }
  \label{fig:sliding_results}
\end{figure}

In Fig. \ref{fig:sliding_spag} we see that over an assimilation of 200 seconds, the ensemble converges on the hidden, true state.
To test the performance of flow reversal prediction, we take the average of the ensemble flow direction (the average of each value of $\phi$) as the predicted flow direction, and count how often we predict reversals both when they do and do not occur.
Varying both the number of model variables and the strength of covariance shifting in Fig. \ref{fig:sliding_results}, we find that covariance shifting improves flow reversal prediction skill even when spatial observation density is decreased.
With full observations (spacing of 1), we obtain a the best predictions with a covariance shift of 2.
For 1/2 and 1/5 observations [a spacing of 2 (5) to observe every other (fifth) variable], we again have the best predictions with a shift of 2.
And for a spacing of 10, observing every 10th variable, we achieve greater prediction skill with a covariance shift of 10.

Computing the average flow direction inside a localized covariance zone is straightforward, and computationally easy since the velocity is immediately available, making incorporation of this scheme into any data assimilation method easy.
Since observations are also sparse in large weather models, we expect that using an adaptive local covariance scheme could lead to improved prediction skill with sparse observations \cite{bishop2011a}.

\subsection*{Dynamic mode decomposition}
\label{dmd_section}

To incorporate limited observations into a high-dimensional CFD simulation, we combine ideas from both CFD literature and data assimilation to make predictions.
We proceed with Tu's algorithm using snapshots every 10 seconds for the first 900 seconds of model time.
A full picture of this time series can be found in \nameref{fig:DMD-timeseries}.

In this reduced space, we extract the modes that correspond to the instability leading to flow reversals.
With a known low-dimensional model of the thermosyphon dynamics, we take this opportunity to test whether DMD can discover the underlying system.
The time series of the model state projection onto a specific DMD mode will represent the time dynamics of a mode that is representative of a single low dimensional variable.

\begin{figure}[h]
  \centering
  \includegraphics[width=0.45\textwidth]{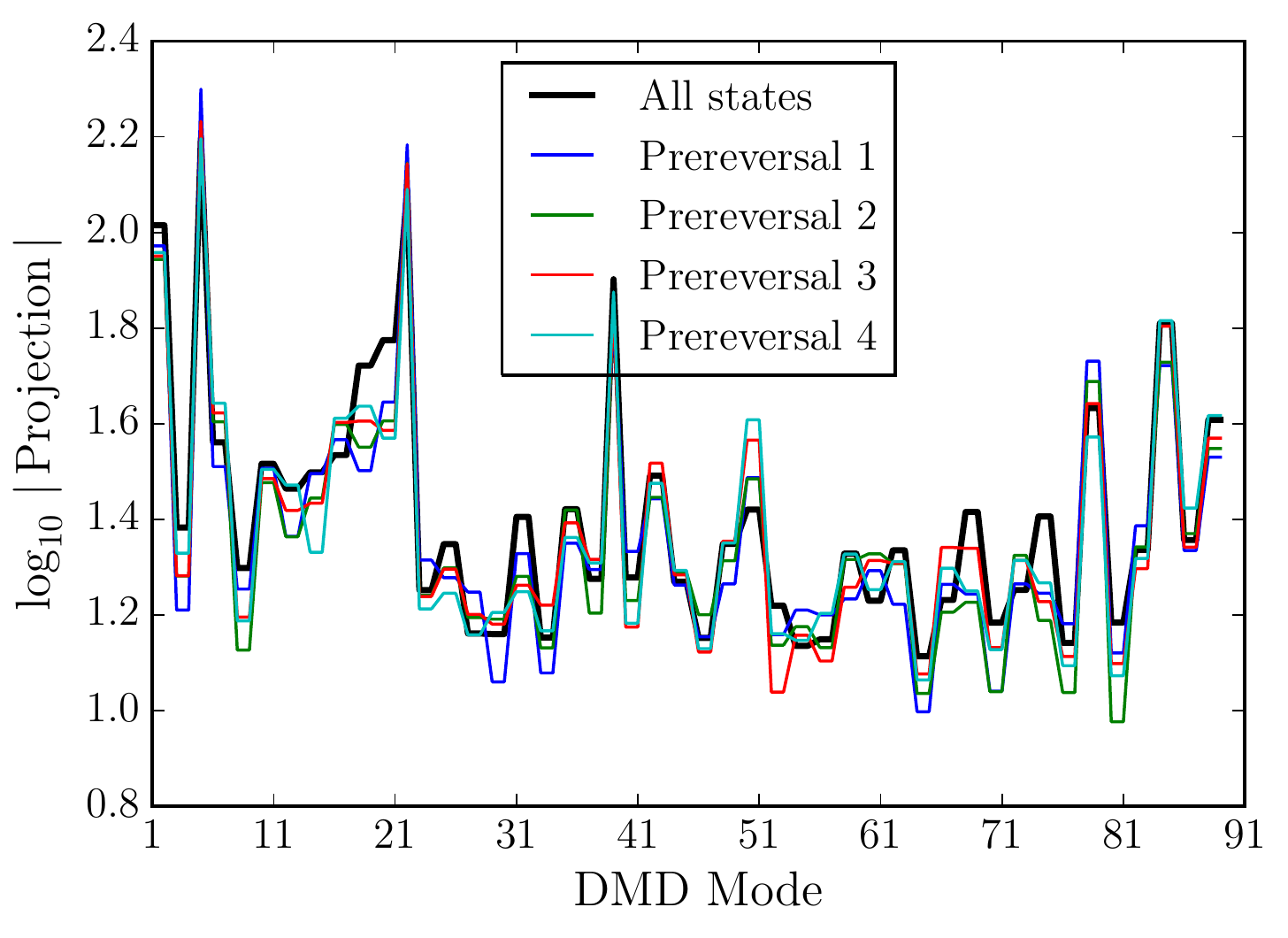}
  \caption[]{
    The $\log_{10}$ average projection onto each DMD mode for different sets of model states.
    DMD constructed as snapshots every 10 seconds for the first 900 seconds of model time, and model states from the first 2000 seconds are all projected onto the DMD modes.
    All states average shown in black, and the average of the subset of states that occur 1 second, 3 seconds, 5 seconds, and 7 seconds before a reversal are shown in other colors.
    The symmetry of the loop generates modes that often come in pairs.
      }
  \label{fig:DMD_modes}
\end{figure}

To look at all of the modes at once, we examine the average magnitude of the projection from all model states onto each mode in comparison to the projection of all states that 1,3,5, and 7 time steps before a reversal.
The magnitude of the mode projection of a predictive mode before a reversal should stand out against the projection average across all states, and decay back towards the average further from the reversal in time.
For modes 21 and 79, we directly observe in Fig. \ref{fig:DMD_modes} that the average projection from states just 1 second before reversal is the most different from the average state projection, and the further away from the reversal the more similar the states become to the average.

\begin{figure}[h]
  \centering
  \includegraphics[width=0.98\textwidth]{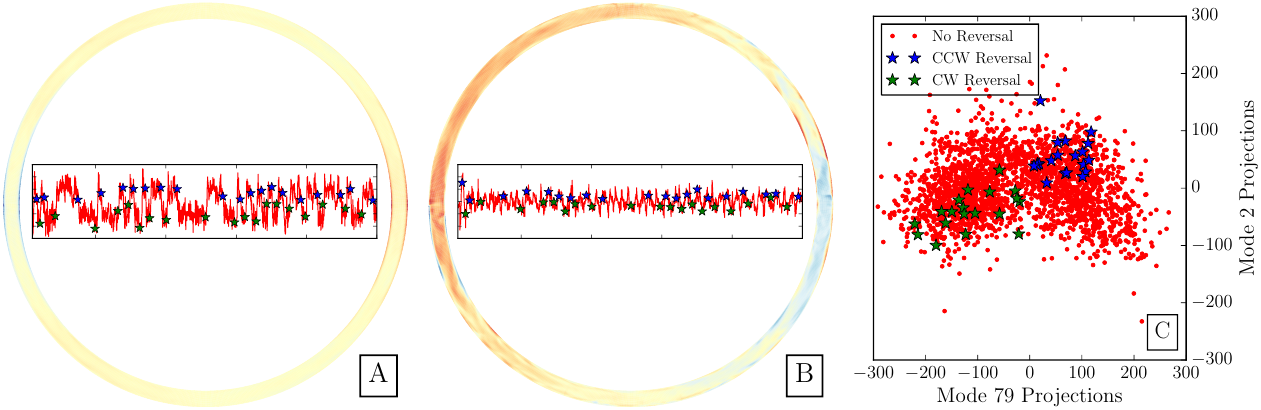}
  \caption[]{
    Panel A: The temperature profile of the thermosyphon of Mode 2, with inset of the projection of time series states onto Mode 2 (the projection coefficient).
    The color scale on the thermosyphon spans the values 1 to 0 in the DMD mode.
    The inset figure is the projection coefficient from time 100 to time 5000, with the projection range being shown from -300 to 300 (as in Panel C) and the starred reversals labeled as in Panel C.
    Panel B: Likewise, the temperature profile of the thermosyphon of Mode 79, with inset of the projection of time series states onto Mode 79 (the projection coefficient).
    The color scale and inset figure axes are the same as Panel A.
    Panel C: A butterfly-shaped phase plane shows the value of the projection onto modes 2 and 79 for each time in the first 2000 time steps of our ground truth model run.
    In blue and green stars the states that occur directly before a flow are highlighted, and are isolated into separate quadrants of phase space.
  }
  \label{fig:DMD_phaseplane}
\end{figure}

We are particularly interested in whether the mode projection time series is predictive of flow reversals, as is true with the hidden system.
The insets of Panel A and Panel B in Fig. \ref{fig:DMD_phaseplane} show two such timeseries, with stars indicating the time of flow reversals.
Individually, these modes increase in amplitude when flow reversals happen.
As a dominant mode, the time series of Mode 2 tracks closely to the timeseries from which the modes were generated, while the dynamics of the projection of Mode 79 are less obvious.

By combining the state projection onto specific mode time series into a phase plane, the combined signal from two modes is used for discovering states that separate reversals in direction and from other states in the phase plane.
In Fig. \ref{fig:DMD_phaseplane} we see that the dominant dynamics from mode 2 plotted with those of mode 79 are able to strongly separate reversals into quadrants of the low-dimensional space.
This result indicates that DMD could be used to improve predictability of reversals.

\section*{Concluding Remarks}

The first output of our work is a general data assimilation framework for MATLAB and Julia.
By utilizing an object-oriented (OO) design, the model and data assimilation algorithm code are separate and can be changed independently.
The principal advantage of this approach is the ease of incorporation of new models and DA techniques (code available at \url{https://github.com/andyreagan/julia-openfoam}).

We next present the results pertaining to the accuracy of forecasts for synthetic data (twin model experiments).
There are many possible experiments given the choice of assimilation window, data assimilation algorithm, localization scheme, model resolution, observational density, observed variables, and observation quality.
We focused on considering the effect of observations and observational locations on the resulting forecast skill, and we find that there is a threshold for the required number of observations to make useful predictions.
In general, and unsurprisingly, we see that increasing observational density leads to improved forecast accuracy.
With too few observations, the data assimilation is unable to recover the underlying dynamics.
Using adaptively localized covariance holds promise for data assimilation with data-scarce models, to overcome the lack of data.

The ability of DMD to recover the lower dimensional dynamics is expected but with 240,000 variables is nonetheless an accomplishment.
When modeling systems for which there are unknown but useful dimension reductions, as demonstrated here, DMD can be a useful tool to find such dimension reductions.
When computational model runs are exceedingly costly or time consuming, the best-guess state projection onto DMD modes provides insights into the system dynamics that could not otherwise be obtained.

The numerical coupling of CFD to experiment by DA should be generally useful to improve the skill of CFD predictions of experiments.
In addition, the CFD model can provide better knowledge of unobservable quantities of interest in fluid flow that use the experimental data to find the analysis state provided by DA.
Adaptive covariance localization further enhances the benefit provided by DA in this context.

\begin{acknowledgements}
    This work was made possible by funding from the Vermont Space Grant Consortium, NASA EPSCoR, NSF-DMS Grant No. 0940271, the Mathematics \& Climate Research Network and the Vermont Advanced Computing Center.
\end{acknowledgements}

\bibliographystyle{unsrt}

\clearpage
\pagebreak

\newwrite\tempfile
\immediate\openout\tempfile=startsupp.txt
\immediate\write\tempfile{\thepage}
\immediate\closeout\tempfile

\onecolumngrid

\setcounter{page}{1}
\renewcommand{\thepage}{S\arabic{page}}
\renewcommand{\thefigure}{S\arabic{figure}}
\renewcommand{\thetable}{S\arabic{table}}
\setcounter{figure}{0}
\setcounter{table}{0}

\appendix

\subsection*{S1 Computational Details and Explicit Equations Used}
\label{S1}

In this section, we first present the governing equations for the flow in our thermal convection loop experiment.
A spatial and temporal discretization of the governing equations is then necessary so that they may be solved numerically.
After discretization, we must specify the boundary conditions.
With the mesh and boundary conditions in place, we can then simulate the flow with a computational fluid dynamics solver.

We now discuss the equations, mesh, boundary conditions, and solver in more detail.
With these considerations, we present our simulations of the thermosyphon.
For a complete derivation of the equations used, see \cite{reagan2013}.

We consider the incompressible Navier-Stokes equations with the Boussinesq approximation to model the flow of water inside a thermal convection loop.
Here we present the main equations that are solved numerically, noting the assumptions that are necessary in their derivation.
In standard notation, for $u,v,w$ the velocity in the $x,y,z$ direction, respectively, the continuity equation for an incompressible fluid is
\begin{equation} \frac{\partial u}{\partial x} + \frac{\partial v}{\partial y} + \frac{\partial w}{\partial z} = 0. \label{eq:NScontIco} \end{equation}

The momentum equations, in tensor notation with bars representing averaged quantities (long timesteps are used, requiring integration), are
\begin{equation} \rho_\text{ref} \left ( \frac{\partial \bar{u}_i}{\partial t} + \frac{\partial}{\partial x_j} \left( \bar{u}_j \bar{u}_i \right) \right )
= -\frac{\partial \bar{p}} {\partial{x_i}} +  \mu \frac{\partial \bar{u}_i}{\partial x_j^2} + \rhoref \left(1 - \beta (T - T_\text{ref})\right) g_i \end{equation}
for $\rho_\text{ref}$ the reference density with the Boussinesq approximation included, $p$ the pressure, $\mu$ the viscosity, and $g_i$ gravity in the $i$-direction.
Note that $g_i = 0$ for $i \in \{ x,y\}$ since gravity is assumed to be the $z$ direction.
Since the model is incompressible, of course our energy equation includes only temperature, and is given by
                            \begin{equation} \partialdiff{T}{t} + \partialdiff{}{x_j} \left ( \rhoref T \overline{u}_j \right ) - \frac{\partial ^2 \alpha T}{\partial x_j \partial x_i}
  =
  - \frac{\partial q_k^*}{\partial x_k}
  - \frac{\partial \overline{q}_k}{\partial x_k}
\end{equation}
for $T$ the temperature and $q$ the flux (where $q = \overline{q} + q^*$ is the averaging notation).

The PISO (Pressure-Implicit with Splitting of Operators) algorithm derives from the work of \cite{issa1986solution}, and is complementary to the SIMPLE (Semi-Implicit Method for Pressure-Linked Equations) \cite{patankar1972calculation} iterative method.
The main difference of the PISO and SIMPLE algorithms is that in the PISO, no under-relaxation is applied and the momentum corrector step is performed more than once \cite{ferziger1996computational}.
They sum up the algorithm in nine steps:
\begin{itemize}
\item Set the boundary conditions.
\item Solve the discretized momentum equation to compute an intermediate velocity field.
\item Compute the mass fluxes at the cell faces.
\item Solve the pressure equation.
\item Correct the mass fluxes at the cell faces.
\item Correct the velocity with respect to the new pressure field.
\item Update the boundary conditions.
\item Repeat from step \#3 for the prescribed number of times.
\item Repeat (with increased time step).
\end{itemize}

The solver itself has 647 dependencies, of which I present only a fraction.
The main code is straight forward, relying on include statements to load the libraries and equations to be solved.
\lstset{language=C++,
        basicstyle=\ttfamily\scriptsize\singlespacing,
        keywordstyle=\color{blue},
        stringstyle=\color{red},
        commentstyle=\color{green},
        morecomment=[l][\color{magenta}]{\#},
        frame=L,
        xleftmargin=\parindent,
                                numbersep=5pt,
        breaklines=true,                breakatwhitespace=false,            escapeinside={\%*}{*)} 
}

\begin{lstlisting}[language=C++,firstline=48,lastline=53]
#include "fvCFD.H"
#include "singlePhaseTransportModel.H"
#include "RASModel.H" // AJR edited header 2013-10-14
#include "radiationModel.H" // (model loaded but not used)
#include "fvIOoptionList.H" // (loaded, but also not used)
#include "pimpleControl.H"
\end{lstlisting}

The main function is then

\begin{lstlisting}[language=C++,firstline=57,lastline=70]
int main(int argc, char *argv[])
{
    #include "setRootCase.H"
    #include "createTime.H"
    #include "createMesh.H"
    #include "readGravitationalAcceleration.H"
    #include "createFields.H"
    #include "createIncompressibleRadiationModel.H"
    #include "createFvOptions.H"
    #include "initContinuityErrs.H"
    #include "readTimeControls.H"
    #include "CourantNo.H"
    #include "setInitialDeltaT.H"
    pimpleControl pimple(mesh);
\end{lstlisting}

We then enter the main loop.
This is computed for each time step, prescribed before the solver is applied.
Note that the capacity is available for adaptive time steps, choosing to keep the Courant number below some threshold, but I do not use this.
For the distributed ensemble of model runs, it is important that each model complete in nearly the same time, so that the analysis is not waiting on one model and therefore under-utilizing the available resources.

\begin{lstlisting}[language=C++,firstline=71,lastline=88]
    while (runTime.loop())
    {
        #include "readTimeControls.H"
        #include "CourantNo.H"
        #include "setDeltaT.H"
        while (pimple.loop())
        {
            #include "UEqn.H"
            #include "TEqn.H"
            while (pimple.correct())
            {
                #include "pEqn.H"
            }
        }
	if (pimple.turbCorr())
	{
	    turbulence->correct();
	}
\end{lstlisting}

Opening up the equation for $U$ we see that Equation

\begin{lstlisting}[language=C++]
    // Solve the momentum equation
    fvVectorMatrix UEqn
    (
        fvm::ddt(U)
      + fvm::div(phi, U)
      + turbulence->divDevReff(U)
     ==
        fvOptions(U)
    );
    UEqn.relax();
    fvOptions.constrain(UEqn);
    if (pimple.momentumPredictor())
    {
        solve
        (
            UEqn
         ==
            fvc::reconstruct
            (
                (
                  - ghf*fvc::snGrad(rhok)
                  - fvc::snGrad(p_rgh)
                )*mesh.magSf()
            )
        );
        fvOptions.correct(U);
    }
\end{lstlisting}

Solving for $T$ is 

\begin{lstlisting}[language=C++]
{
    // with our laminar model, alphat = 0
    alphat = turbulence->nut()/Prt;
    alphat.correctBoundaryConditions();
    volScalarField alphaEff("alphaEff", turbulence->nu()/Pr + alphat);
    fvScalarMatrix TEqn
    (
        fvm::ddt(T)
      + fvm::div(phi, T)
      - fvm::laplacian(alphaEff, T)
     ==
        radiation->ST(rhoCpRef, T)
      + fvOptions(T)
    );
    TEqn.relax();
    fvOptions.constrain(TEqn);
    TEqn.solve();
    radiation->correct();
    fvOptions.correct(T);
    rhok = 1.0 - beta*(T - TRef); // Boussinesq approximation
}
\end{lstlisting}

Finally, we solve for the pressure $p$ in ``pEqn.H'':

\begin{lstlisting}[language=C++]
{
    volScalarField rAU("rAU", 1.0/UEqn.A());
    surfaceScalarField rAUf("Dp", fvc::interpolate(rAU));
    volVectorField HbyA("HbyA", U);
    HbyA = rAU*UEqn.H();
    surfaceScalarField phig(-rAUf*ghf*fvc::snGrad(rhok)*mesh.magSf());
    surfaceScalarField phiHbyA
    (
        "phiHbyA",
        (fvc::interpolate(HbyA) & mesh.Sf())
      + fvc::ddtPhiCorr(rAU, U, phi)
      + phig
    );
    while (pimple.correctNonOrthogonal())
    {
        fvScalarMatrix p_rghEqn
        (
            fvm::laplacian(rAUf, p_rgh) == fvc::div(phiHbyA)
        );
        p_rghEqn.setReference(pRefCell, getRefCellValue(p_rgh, pRefCell));
        p_rghEqn.solve(mesh.solver(p_rgh.select(pimple.finalInnerIter())));
        if (pimple.finalNonOrthogonalIter())
        {
            // Calculate the conservative fluxes
            phi = phiHbyA - p_rghEqn.flux();
            // Explicitly relax pressure for momentum corrector
            p_rgh.relax();
            // Correct the momentum source with the pressure gradient flux
            // calculated from the relaxed pressure
            U = HbyA + rAU*fvc::reconstruct((phig - p_rghEqn.flux())/rAUf);
            U.correctBoundaryConditions();
        }
    }
    #include "continuityErrs.H"
    p = p_rgh + rhok*gh;
    if (p_rgh.needReference())
    {
        p += dimensionedScalar
        (
            "p",
            p.dimensions(),
            pRefValue - getRefCellValue(p, pRefCell)
        );
        p_rgh = p - rhok*gh;
    }
}
\end{lstlisting}

The final operation being the conversion of pressure to hydrostatic pressure,
\begin{equation*} p _\text{rgh} = p - \rho _k g_h . \end{equation*}
This ``pEqn.H'' is then re-run until convergence is achieved, and the PISO loop begins again.

\begin{figure}[h!]
  \centering
    \includegraphics[width=0.79\textwidth]{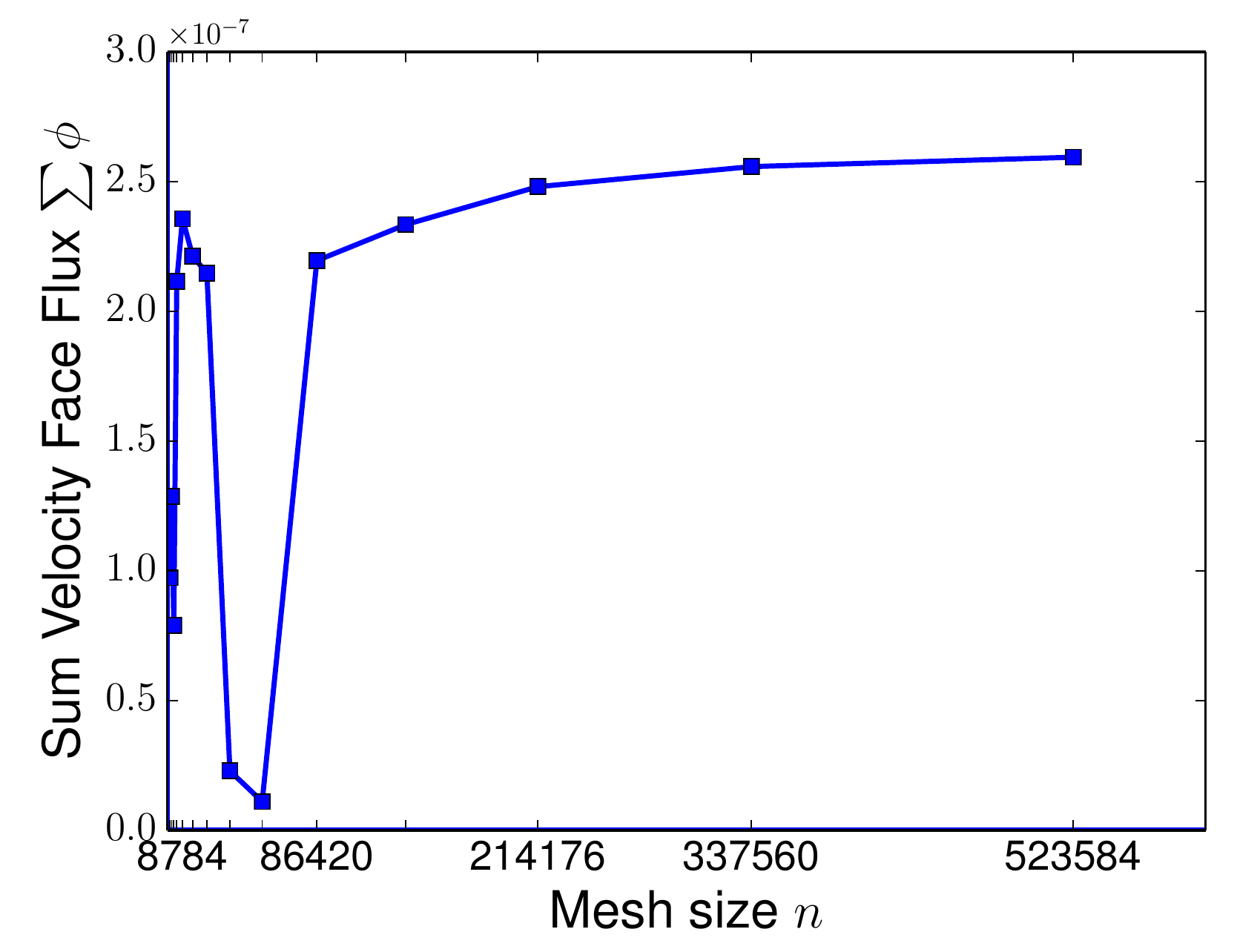}
  \caption[Long-term behavior for different meshes]
  {
    With a fixed choice of solver, boundary conditions, and initial conditions that lead to a stable convective state, we present the long-term behavior of the velocity face flux at the top slice for different meshes.
    The face flux is reported as the average for the last 20 times saves for which the velocity flux is summed across a slice perpendicular to the loop, here we show the top slice.
    We choose a fixed time step of 0.005 for each simulation, and run the solver for 60 hours on 8 cores.
    The computational limit of mesh creation was a memory limit at 818280 cells, so we present results for meshes starting at 1600 cells and cells decreasing in size by a factor of 1.25 in both $y$ and $z$ up to a mesh of 523584 cells.
    For meshes with more than 80,000 cells we see that the solutions are very similar.
    The smaller meshes generate increasing unstable flow behavior, leading to oscillations of flux and then flow reversals for the smallest meshes of size 2500 and 1600 cells.
  }
  \label{fig:meshverification}
\end{figure}

We verify convergence of the solution in a steady flow state regime in Fig \ref{fig:meshverification}.

\clearpage
\pagebreak
\subsection*{S2 The Ehrhard and M\"{u}ller Equations}
\label{S2}

Following the derivation by Harris \cite{harris2011predicting}, itself a representation of the derivation of Gorman \cite{gorman1986} and namesakes Ehrhard and M\"{u}ller \cite{ehrhard1990dynamical}, we derive the equations governing a closed loop thermosyphon.

Similar to the derivation of the governing equations of computational fluid dynamics, we start with a small but finite volume inside the loop.
Here, however, the volume is described by $\pi r^2 R \text{d} \phi$ for $r$ the interior loop size (such that $\pi r^2$ is the area of a slice) and $R\text{d}\phi$ the arc length (width) of the slice.
Newton's second law states that momentum is conserved, such that the sum of the forces acting upon our finite volume is equal to the change in momentum of this volume.
Therefore we have the basic starting point for forces $\sum F$ and velocity $u$ as
\begin{equation} \sum F = \rho \pi r^2 R \text{d}\phi \diff{u}{t} .\end{equation}
The sum of the forces is $\sum F = F_{\{p,f,g\}}$ for net pressure, fluid shear, and gravity, respectively.
We write these as
\begin{align} & F_p = -\pi r^2 \text{d} \phi \pdiff{p}{\phi}\\
& F_w = -\rho \pi r^2 \text{d} \phi f_w\\
& F_g = -\rho \pi r^2 \text{d} \phi g \sin (\phi)\end{align}
where $\partial p /\partial \phi$ is the pressure gradient, $f_w$ is the wall friction force, and $g \sin (\phi)$ is the vertical component of gravity acting on the volume.

We now introduce the Boussinesq approximation which states that both variations in fluid density are linear in temperature $T$ and density variation is insignificant except when multiplied by gravity.
The consideration manifests as
\begin{equation*} \rho = \rho (T) \simeq \rho _\text{ref} (1 - \beta (T - T_\text{ref}) \end{equation*}
where $\rho _0$ is the reference density and $T_\text{ref}$ is the reference temperature, and $\beta$ is the thermal expansion coefficient.
The second consideration of the Boussinesq approximation allows us to replace $\rho$ with this $\rhoref$ in all terms except for $F_g$.
We now write momentum equation as
\begin{equation} -\pi r^2 \dphi \pdiff{p}{\phi} - \rhoref \phi r^2 R \dphi f_w
- \rhoref (1 - \rho (T- T_\text{ref}) ) \pi r^2 R \dphi g \sin (\phi) = \rhoref \pi r^2 R \dphi \diff{u}{t}. \end{equation}
Canceling the common $\pi r^2$, dividing by $R$, and pulling out $\dphi$ on the LHS we have
\begin{equation} -\dphi \left ( \pdiff{p}{\phi}  \frac{1}{R} - \rhoref f_w - \rhoref (1 - \rho (T- T_\text{ref}) ) g \sin (\phi) \right ) = \rhoref \dphi \diff{u}{t}. \label{eq:EM07} \end{equation}
We integrate this equation over $\phi$ to eliminate many of the terms, specifically we have
\begin{align*}
& \int _{0} ^{2\pi} -\dphi \pdiff{p}{\phi} \frac{1}{R} \rightarrow 0\\
& \int _{0} ^{2\pi} -\dphi \rhoref g \sin (\phi) \rightarrow 0\\
& \int _{0} ^{2\pi} -\dphi \rhoref \beta T_\text{ref} g \sin (\phi) \rightarrow 0.\end{align*}
Since $u$ (and hence $\diff{u}{\phi}$) and $f_w$ do not depend on $\phi$, we can pull these outside an integral over $\phi$ and therefore the momentum equation is now 
\begin{equation*} 2\pi f_w \rho _0 + \int _{0} ^{2\pi} \dphi \rhoref \beta T g \sin (\phi) = 2\pi \diff{u}{\phi} \rhoref .\end{equation*}
Diving out $2\pi$ and pull constants out of the integral we have our final form of the momentum equation
\begin{equation} f_w \rhoref + \frac{\rhoref \beta g }{2 \pi} \int _{0} ^{2\pi} \dphi T \sin (\phi) = \diff{u}{\phi} \rhoref \label{eq:EM10}.\end{equation}
Now considering the conservation of energy within the thermosyphon, the energy change within a finite volume must be balanced by transfer within the thermosyphon and to the walls.
The internal energy change is given by
\begin{equation} \rhoref \pi r^2 R \dphi \left ( \pdiff{T}{t} + \frac{u}{R}\pdiff{T}{\phi} \right ) \label{eq:EMeg1}\end{equation}
which must equal the energy transfer through the wall, which is, for $T_w$ the wall temperature:
\begin{equation} \dot{q} = -\pi r^2 R \dphi h_w (T - T_w) . \label{eq:EMeg2} \end{equation}
Combining Equations \ref{eq:EMeg1} and \ref{eq:EMeg2} (and canceling terms) we have the energy equation:
\begin{equation} \left ( \pdiff{T}{t} + \frac{u}{R}\pdiff{T}{\phi} \right ) = \frac{-h_w}{\rhoref c_p} \left( T - T_w \right ) \label{eq:EMeq}.\end{equation}
The $f_w$ which we have yet to define and $h_w$ are fluid-wall coefficients and can be described by \cite{ehrhard1990dynamical}:
\begin{align*} & h_w = h_{w_0} \left ( 1 + K h(|x_1|) \right ) \\
& f_w = \frac{1}{2} \rhoref f_{w_0} u .\end{align*}
We have introduced an additional function $h$ to describe the behavior of the dimensionless velocity $x_1 \alpha u$.
This function is defined piece-wise as 
\begin{equation} h (x) = \left \{ \begin{array}{ll} x^{1/3} & ~~\text{when} ~x \geq 1\\ p (x) & ~~\text{when} ~ x <1 \end{array} \right. \label{eq:h_defined} \end{equation} 
where $p(x)$ can be defined as $p(x) = \left( 44x^2 -55 x^3 + 20x^4 \right ) /9$ such that $p$ is analytic at 0 \cite{harris2011predicting}.

Taking the lowest modes of a Fourier expansion for $T$ for an approximate solution, we consider:
\begin{equation} T(\phi , t) = C_0 (t) + S(t) \sin (\phi ) + C(t) \cos (\phi) . \end{equation}
By substituting this form into Equations \ref{eq:EM10} and \ref{eq:EMeq} and integrating, we obtain a system of three equations for our solution.
We then follow the particular nondimensionalization choice of Harris \etal such that we obtain the following ODE system, which we refer to as the Ehrhard-M\"{u}ller equations:
\begin{align}
& \diff{x_1}{t'} = \alpha (x_2 - x_1),\\
& \diff{x_2}{t'} = \beta x_1 - x_2 (1 + Kh(|x_1|)) - x_1x_3,\\
& \diff{x_3}{t'} = x_1x_2 - x_3 (1 + Kh(|x_1|)) .\end{align}
The nondimensionalization is given by the change of variables
\begin{align}
& t' = \frac{h_{w_0}}{\rhoref c_p}t,\\
& x_1 = \frac{\rhoref c_p }{R h_{w_0}} u, \\
& x_2 = \frac{1}{2} \frac{\rhoref c_p \beta g}{ R h_{w_0} f_{w_0}} \Delta T_{3-9}, \\
& x_3 = \frac{1}{2} \frac{\rhoref c_p \beta g}{ R h_{w_0} f_{w_0}} \left ( \frac{4}{\pi} \Delta T_w - \Delta T_{6-12} \right ) 
\end{align}
and
\begin{align}
& \alpha = \frac{1}{2} R c_p f_{w_0} / h_{w_0} ,\\
& \gamma = \frac{2}{\pi} \frac{\rhoref c_p \beta g}{Rh_{w_0} f_{w_0}} \Delta T_w. \end{align}

Through careful consideration of these non-dimensional variable transformations we verify that $x_1$ is representative of the mean fluid velocity, $x_2$ of the temperature difference between the 3 and 9 o'clock positions on the thermosyphon, and $x_3$ the deviation from the vertical temperature profile in a conduction state \cite{harris2011predicting}.

\clearpage
\pagebreak
\subsection*{S3 Data Assimiliation}
\label{S3}

The TLM is the model which advances an initial perturbation $\delta \mbx_{i}$ at timestep $i$ to a final perturbation $\delta \mbx_{i+1}$ at timestep $i+1$.
The dynamical system we are interested in, Lorenz '63, is given as a system of ODE's:
\[ \frac{d\mbx}{dt} = F(\mbx) .\]
We integrate this system using a numerical scheme of our choice (in the given examples we use a second-order Runge-Kutta method), to obtain a model $M$ discretized in time.
\[ \mbx(t) = M[ \mbx(t_0) ] .\]
Introducing a small perturbation $\mby$, we can approximate our model $M$ applied to $\mbx(t_0) + \mby(t_0)$ with a Taylor series around $\mbx(t_0)$:
\begin{align*} M[ \mbx(t_0) + \mby(t_0) ] &= M [ \mbx(t_0) ] + \frac{\partial M}{\partial \mbx} \mby(t_0) + O [ \mby(t_0) ^2 ]\\ &\approx \mbx(t) + \frac{\partial M}{\partial \mbx} \mby(t_0) .\end{align*}
We can then solve for the linear evolution of the small perturbation $\mby(t_0)$ as 
\begin{equation} \frac{d\mby }{dt } = \mathbf{J} \mby \label{eq:ODETLM} \end{equation}
where $\mathbf{J} = \partial F / \partial \mbx$ is the Jacobian of $F$.
We can solve the above system of linear ordinary differential equations using the same numerical scheme as we did for the nonlinear model.

\begin{figure}[h]
  \centering
  \includegraphics[width=0.89\textwidth]{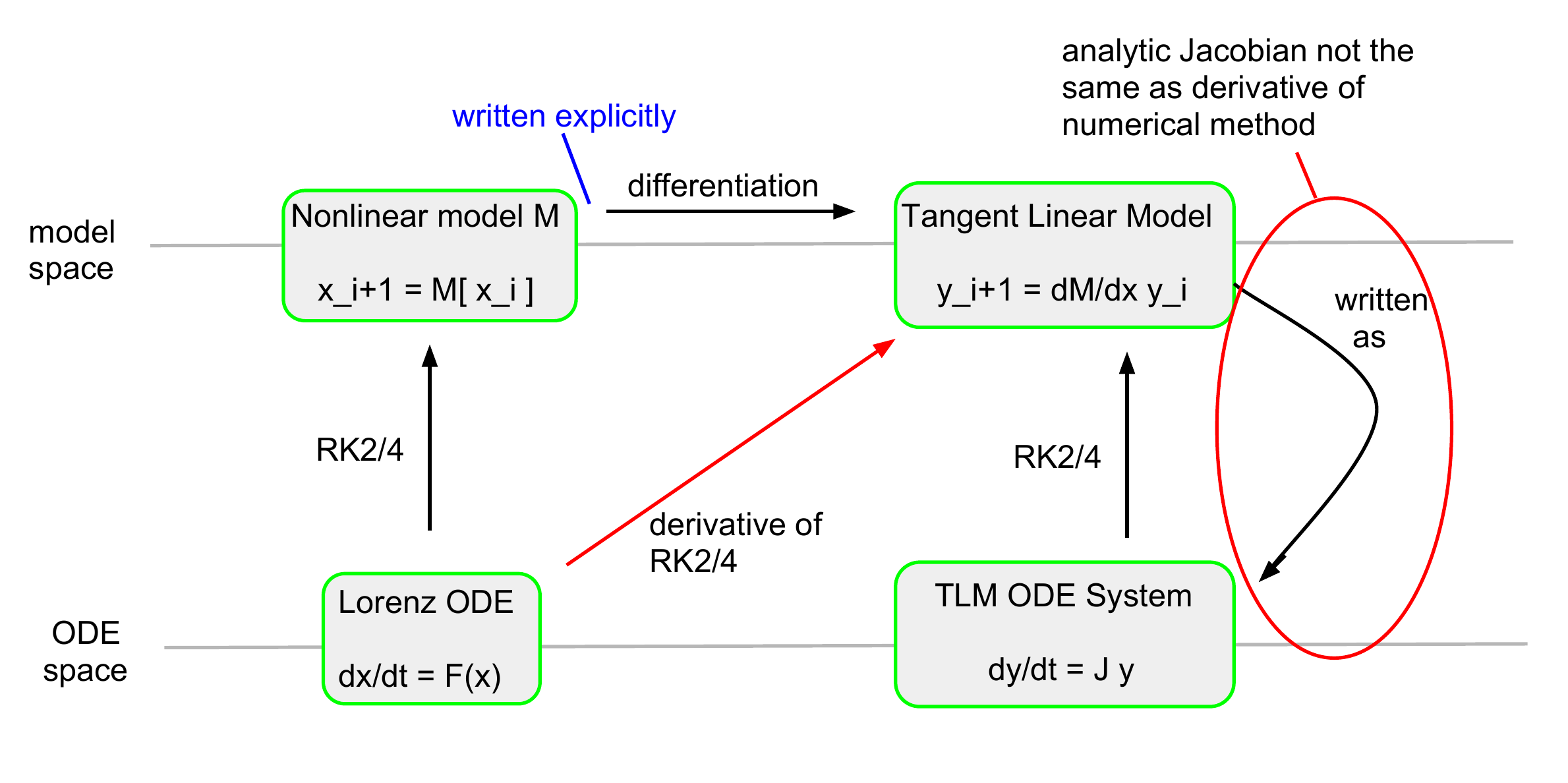}
  \caption[An explanation of how and why the best way to obtain a TLM is with a differentiated numerical scheme]{
    An explanation of how and why the best way to obtain a TLM is with a differentiated numerical scheme.
    Both the Lorenz ODE and TLM ODE System can be solved by RK2/4, but the analytic Jacobian of TLM that would is not the same as the derivative of the numerical method.
    In particular, the derivative of the RK2/4 integrator is used to obtain a TLM that most accurately propogates error growth in the Lorenz '63 system.
  }
  \label{fig:TLMscheme}
\end{figure}

One problem with solving the system of equations given by Equation \ref{eq:ODETLM} is that the Jacobian matrix of discretized code is not necessarily identical to the discretization of the Jacobian operator for the analytic system.
This is a problem because we need to have the TLM of our model $M$, which is the time-space discretization of the solution to $d\mbx/dt = F(\mbx)$.
We can apply our numerical method to the $d\mbx/dt = F(\mbx)$ to obtain $M$ explicitly, and then take the Jacobian of the result.
This method is, however, prohibitively costly, since Runge-Kutta methods are implicit.
It is therefore desirable to take the derivative of the numerical scheme directly, and apply this differentiated numerical scheme to the system of equations $F(\mbx)$ to obtain the TLM.
A schematic of this scenario is illustrated in Figure \ref{fig:TLMscheme}.
To that the derivative of numerical code for implementing the EKF on models larger than 3 dimensions (i.e. global weather models written in Fortan), automatic code differentiation is used \cite{autodiff1981}.

\begin{figure}[h]
  \centering
  \includegraphics[width=0.45\textwidth]{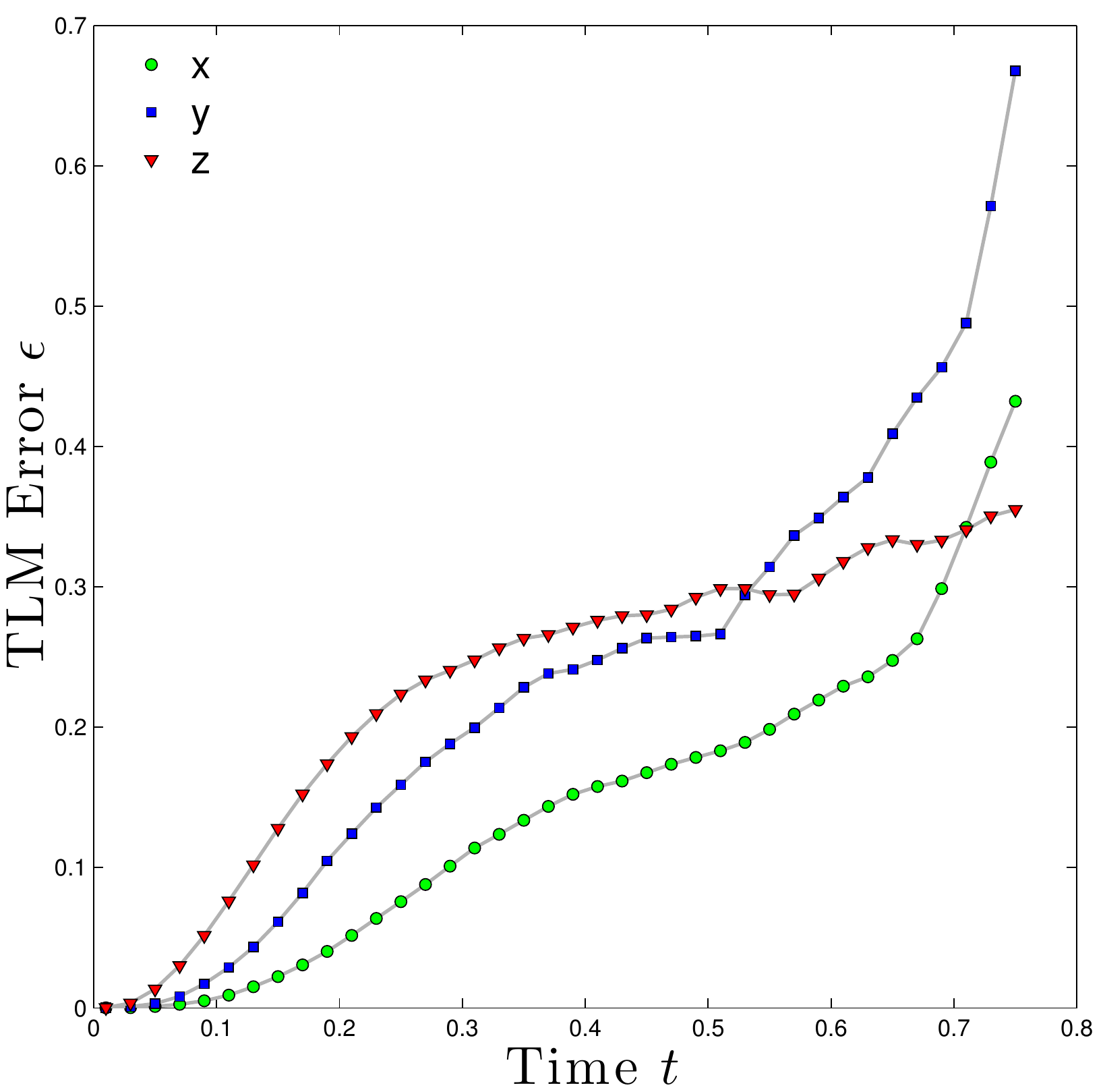}
  \caption[The future error predicted by the TLM is compared to the error growth in Lorenz '63 system for an initial perturbation with standard deviation of 0.1, averaged over 1000 TLM integrations]{
    The future error predicted by the TLM is compared to the error growth in Lorenz '63 system for an initial perturbation with standard deviation of 0.1, averaged over 1000 TLM integrations.
    The $\epsilon$ is not the error predicted by the TLM, but rather the error of the TLM in predicting the error growth.
    We see an intially linear error growth for small time, which is overcome by the nonlinearity of the Lorenz system for longer time.
  }
  \label{fig:TLMverification}
\end{figure}

To verify our implementation of the TLM, we propagate a small error in the Lorenz '63 system and plot the difference between that error and the TLM predicted error, for each variable (Figure \ref{fig:TLMverification}).

\begin{figure}[h]
  \centering
  \includegraphics[width=0.79\textwidth]{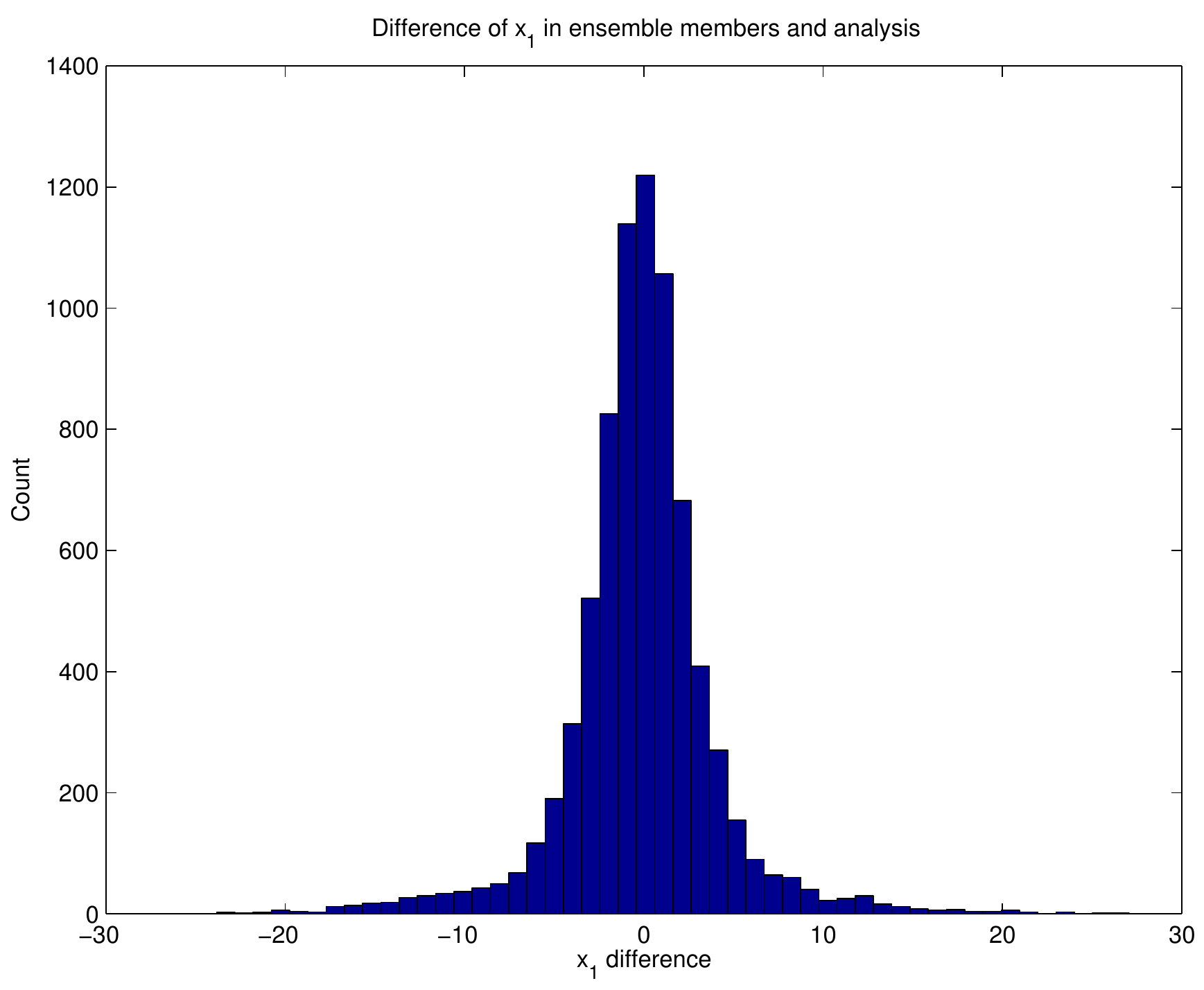}
  \caption[The difference of ensemble forecasts from the analysis is reported for 760 assimilation windows in one model run of length 200, with 10 ensemble members and an assimilation window of length 0.261]{
    The difference of ensemble forecasts from the analysis is reported for 760 assimilation windows in one model run of length 200, with 10 ensemble members and an assimilation window of length 0.261.
    This has the same shape of as the difference between ensemble forecasts and the mean of the forecasts (not shown).
    This spread of ensemble forecasts is what allows us to estimate the error covariance of the forecast model, and appears to be normally distributed.
  }
  \label{fig:EnKFhist}
\end{figure}

With a finite ensemble size, the ensemble method is only an approximation and therefore in practice it often fails to capture the full spread of error.
To better capture the model variance, additive and multiplicative inflation factors are used to obtain a good estimate of the error covariance matrix (\nameref{fig:ETKF_cov_tuning_390s} Section).
The spread of ensemble members in the $x_1$ variable of the Lorenz model, as distance from the analysis, can be seen in Figure \ref{fig:EnKFhist}.

In computing the error covariance $\mbP_f$ from the ensemble, we wish to add up the error covariance of each forecast with respect to the mean forecast. 
But this would underestimate the error covariance, since the forecast we're comparing against is used in the ensemble average (to obtain the mean forecast).
Therefore, to compute the error covariance matrix for each forecast, that forecast itself is excluded from the ensemble average forecast.

\begin{figure}[h]
  \centering
  \includegraphics[width=0.99\textwidth]{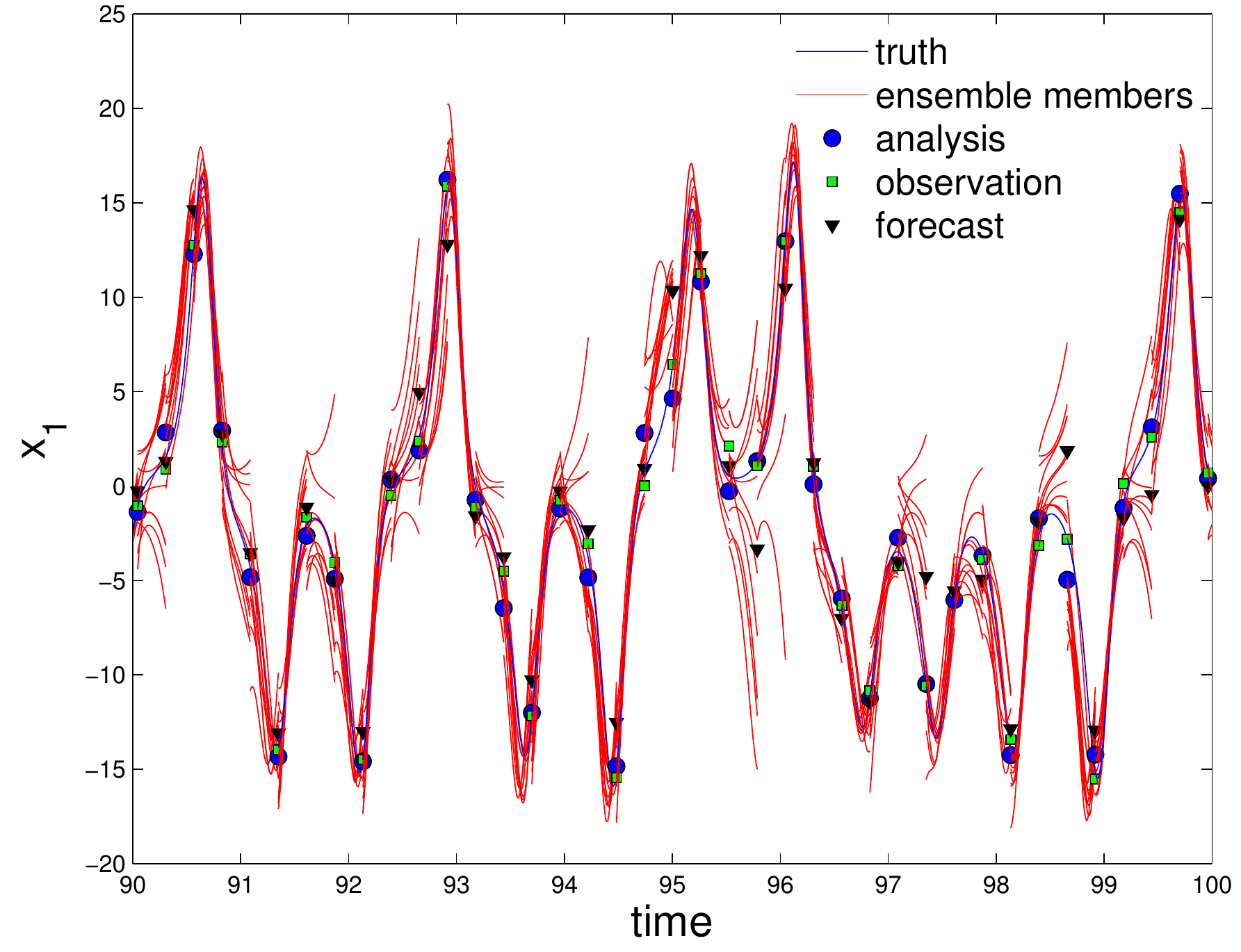}
  \caption[A sample time-series of the ensembles used in the EnKF]{
    A sample time-series of the ensembles used in the EnKF.
    In all tests, as seen here, 10 ensemble members are used.
    For this run, 384 assimilation cycles are performed with a window length of 0.26 model time units.
    We can see that the ensemble member state after assimilation better represents the uncertainty of the analysis state and enables some ensemble members to stay close to the true state.
  }
  \label{fig:spaghetti}
\end{figure}

We can see the classic spaghetti of the ensemble with this filter implemented on Lorenz 63 in Figure \ref{fig:spaghetti}.

We denote the forecast within an ensemble filter as the average of the individual ensemble forecasts, and an explanation for this choice is substantiated by Burgers \cite{burgers1998analysis}.
The general EnKF which we use is most similar to that of Burgers.
Many algorithms based on the EnKF have been proposed and include the Ensemble Transform Kalman Filter (ETKF) \cite{ott2004local}, Ensemble Analysis Filter (EAF) \cite{anderson2001new}, Ensemble Square Root Filter (EnSRF) \cite{tippett2003ensemble}, Local Ensemble Kalman Filter (LEKF) \cite{ott2004local}, and the Local Ensemble Transform Kalman Filter (LETKF) \cite{hunt2007efficient}.
A comprehensive overview through 2003 is provided by Evensen \cite{evensen2003ensemble}.
For further details on the most advanced methods, beyond what is provided in the body of the paper, we direct the reader the above references and the derivations provided in \cite{reagan2013}.

\clearpage
\pagebreak
\subsection*{S4 Additional DMD Details}
\label{S4}

\begin{figure}[h]
  \centering
  \includegraphics[width=0.98\textwidth]{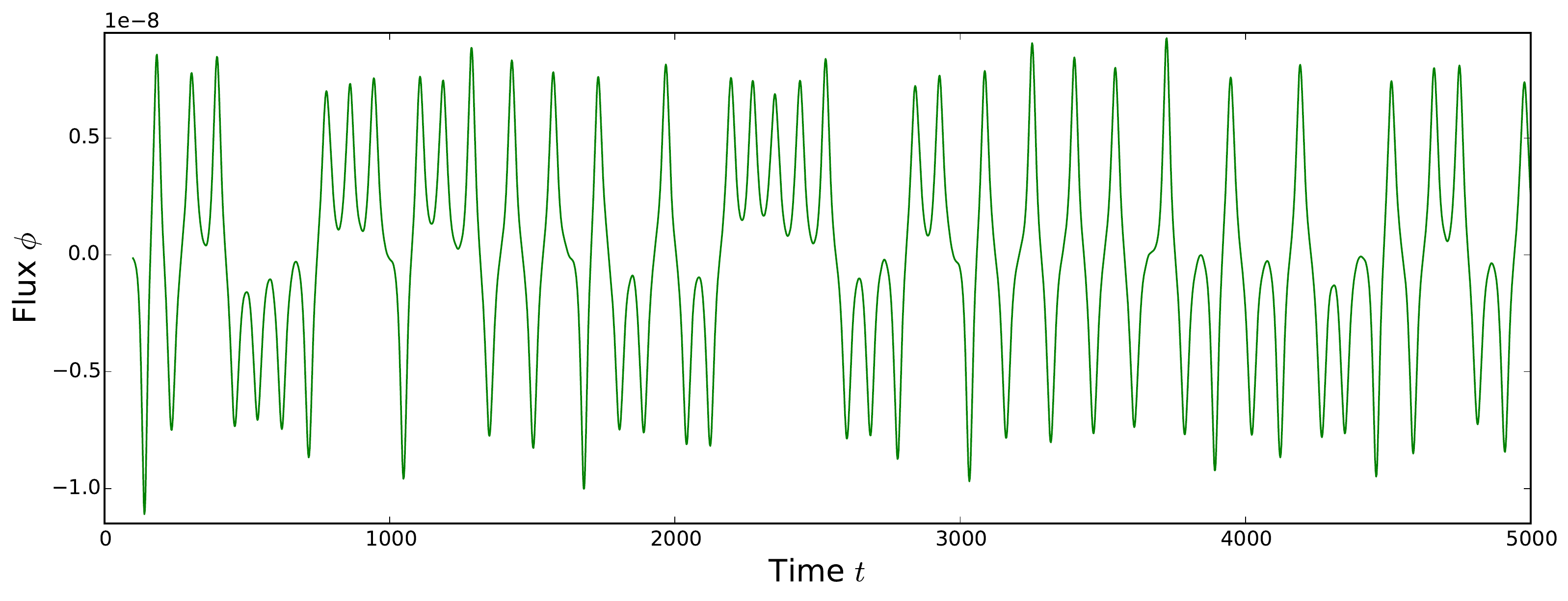}
  \caption[]{
    Flux timeseries on which DMD is performed.
    We report the flux as the sum of the face flux values on a slice of the loop at the 9 o'clock position.
    In this flux timeseries we see dynamics visually similar to the $x_1$ variable of the Lorenz 63 system.
    Residence time in either flow direction is aperiodic and unstable with the flow speed oscillating within a single direction with growing amplitude until reversing.
  }
  \label{fig:DMD-timeseries}  
\end{figure}

The general algorithm for DMD is provided in the \nameref{dmd_section} Section, and here we supply more results of the DMD procedure.
The timeseries from which we computed the decomposition is shown in Figure \ref{fig:DMD-timeseries}.

\begin{figure}[h]
  \centering
  \includegraphics[width=0.98\textwidth]{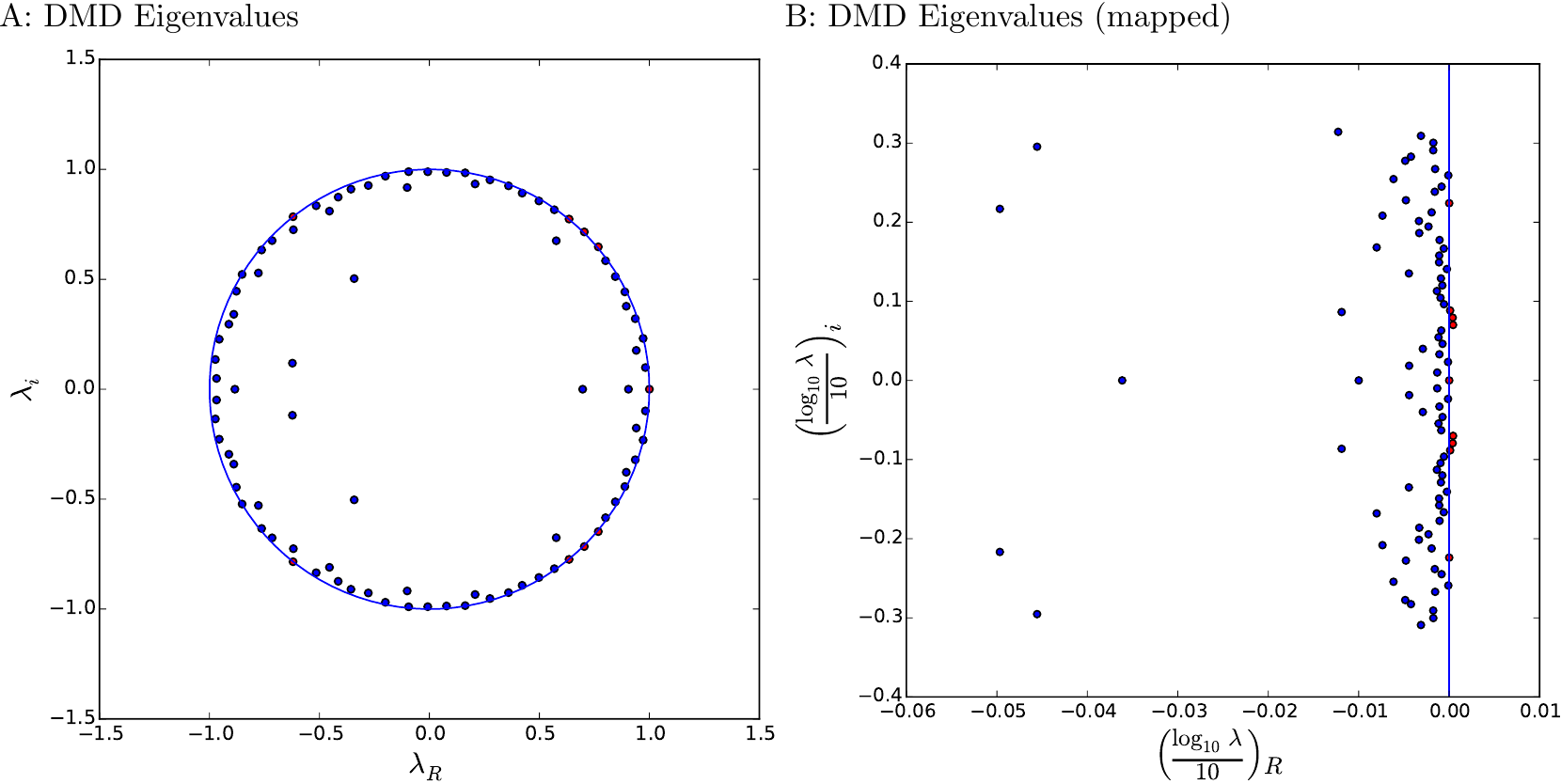}
  \caption[]{
    Eigenvalues of DMD Modes presented in both raw and scaled (mapped) form.
    In Panel A we see the eigenvalue of each DMD mode on the complex plane, with an inset unit circle.
    Those eigenvalues with magnitude greater than 1 are shown in red.
    In Panel B we see the same eigenvalues on the complex plane, transformed by the base 10 logarithm.
    Again we color in red those eigenvalues with real part greater than 0.
  }
  \label{fig:DMD-eigenvalues}  
\end{figure}

The real and imaginary components of the DMD eigenvalues are shown in both un-mapped and mapped form in Figure \ref{fig:DMD-eigenvalues}.

\end{document}